\documentclass[12pt]{article} 
\usepackage{amssymb} 
\usepackage{amsmath} 
\usepackage{latexsym}\usepackage{verbatim} 
\newtheorem{thm}{Theorem}
\newtheorem{prop}[thm]{Proposition}   
\newtheorem{la}[thm]{Lemma} 
   
\newtheorem{dftemp}[thm]{Definition} 
\newtheorem{extemp}[thm]{Example}   
\newtheorem{rmktemp}[thm]{Remark} 
\newtheorem{convtemp}[thm]{Convention} 
\newtheorem{questemp}[thm]{Question} 
\newenvironment{df}{\begin{dftemp}\normalfont}{\end{dftemp}} 
 
\newenvironment{rmk}{\begin{rmktemp}\normalfont}{\end{rmktemp}} 
 
\newenvironment{ques}{\begin{questemp}\normalfont}{\end{questemp}}

\newenvironment{eqs}{\begin{eqnarray*}}{\end{eqnarray*}} 
 
\newenvironment{ls}{\begin{itemize}}{\end{itemize}} 
 
\newenvironment{pf}{\smallskip\noindent\emph{Proof}\quad}%
{\hfill$\square$\par\smallskip} 
 
\newcommand{\ger}[1]{\ensuremath{\mathfrak {#1}}}

\newcommand{\bbb}[1]{\ensuremath{\mathbb {#1}}} 
\newcommand{\ttt}[1]{\ensuremath{\mathtt {#1}}} 
\newcommand{\dom}[1]{\ensuremath{{\text{Dom}}(#1)}} 
\newcommand{\ran}[1]{\ensuremath{{\text{Range}}(#1)}} 
 
\newcommand{\emp}{\varnothing} 
\renewcommand{\phi}{\varphi} 
\newcommand{\eps}{\varepsilon} 
\newcommand{\sq}[1]{\ensuremath{\langle#1\rangle}}

\newcommand{\notarrow}{\kern .42em\not\kern -.42em\longrightarrow} 
\newcommand{\cpt}{\mbox{\~CPT}} 
\newcommand{\cptc}{\ensuremath{\cpt}+\text{Card}} 
\newcommand{\HF}[1]{\ensuremath{\text{HF}(#1)}} 
\newcommand{\od}[1]{\ensuremath{\text{Odd}(#1)}} 
\newcommand{\up}{\ensuremath{\Upsilon}} 
\newenvironment{eatab} 
{\bigskip\noindent\begin{minipage}{\textwidth}\upshape\ttfamily 
 \begin{tabbing}mmm\=mmm\=mmm\=mmm\=mmm\=\kill} 
{\end{tabbing}\end{minipage}\bigskip} 
 
\title{On Polynomial Time Computation\\ 
Over Unordered Structures} 
 
\author{Andreas Blass%
\thanks{Partially supported by NSF grant DMS--0070723 and a grant from 
Microsoft Research.  Address: 
Mathematics Department, University of Michigan, 
Ann Arbor, MI \ 48109--1109, U.S.A{}. email: \texttt{ablass@umich.edu}} 
\and 
Yuri Gurevich%
\thanks{Microsoft Research, 
One Microsoft Way, Redmond, WA \ 98052,
U.S.A{}. e-mail: \texttt{gurevich@microsoft.com}} 
\and 
Saharon Shelah%
\thanks{Partially supported by U.S.-Israel Binational Science
Foundation.  Address:  
Institute of Mathematics, Hebrew University of Jerusalem, 
Givat Ram, 91904 Jerusalem, Israel, 
and Mathematics Department, Rutgers University, 
New Brunswick, NJ \ 08903, U.S.A{}. e-mail:
\texttt{shelah@math.huji.ac.il}}}  
 
\begin{document} 
\maketitle 
 
\begin{abstract} 
This paper is motivated by the question whether there 
exists a logic capturing polynomial time computation over unordered 
structures.  We consider several algorithmic problems near the border 
of the known, logically defined complexity classes contained in 
polynomial time.  We show that fixpoint logic plus counting is stronger 
than might be expected, in that it can express the existence of a 
complete matching in a bipartite graph.  We revisit the known examples 
that separate polynomial time from fixpoint plus counting.  We show 
that the examples in a paper of Cai, F\"urer, and Immerman, when 
suitably padded, are in choiceless polynomial time yet not in fixpoint 
plus counting.  Without padding, they remain in polynomial time but 
appear not to be in choiceless polynomial time plus counting.  Similar 
results hold for the multipede examples of Gurevich and Shelah, except 
that their final version of multipedes is, in a sense, already 
suitably padded.  Finally, we describe another plausible candidate, 
involving determinants, for the task of separating polynomial time 
from choiceless polynomial time plus counting. 
\end{abstract} 
 
\section{Introduction}	\label{intro} 
 
We shall be concerned with computational problems whose inputs are 
finite structures (for a fixed, finite vocabulary \up) and whose 
outputs are ``yes'' and ``no'' (or 1 and 0, or \ttt{true} and 
\ttt{false}).    
 
When \up\ contains a binary relation symbol $\preceq$ interpreted in 
all input structures as a linear ordering of the underlying set, then 
these structures admit an easy, canonical encoding as strings.  In 
this situation, one defines polynomial time computation on ordered 
structures to mean polynomial time Turing machine computation using as 
inputs the string encodings of the structures.  Of course, polynomial 
time is robust, so equivalent definitions could be given using other 
computation models in place of Turing machines.    
 
Turing machines that include a clock to stop the computation after a 
specified polynomial number of steps thus form a computation model 
capturing PTime on ordered structures.  They constitute a logic in the 
broad sense defined in \cite{logic}.  Immerman \cite{imm} and Vardi 
\cite{vardi} showed that PTime on ordered structures is also captured 
by a logic with the look and feel traditionally associated with 
logics, namely fixpoint logic FP.  (We shall review in 
Section~\ref{background} the definitions of FP and other logics 
mentioned in this introduction.) 
 
For unordered input structures, the situation is quite different.  One
can encode such a structure as a string by first choosing a linear
ordering of the underlying set.  Thus, the same structure has many
string encodings, and no efficient way is known to choose a preferred
encoding.  Following Chandra and Harel \cite{c-h}, one says that a
problem having unordered structures as inputs is solvable in
polynomial time if there is a PTime Turing machine that solves the
problem when given any string encoding of the input structure (arising
from any ordering of the underlying set).
 
This does not provide a logic in the sense of \cite{logic}, because 
that sense requires the sentences of a logic to form a recursive set. 
In the case at hand, the ``sentences'' would be PTime Turing machines 
whose output is the same for any two inputs encoding the same 
structure.  This invariance property is undecidable, so the 
recursivity requirement is violated.    
 
Nor does fixpoint logic FP capture PTime on unordered structures.  It 
cannot even express ``the universe has an even number of elements'' 
when the vocabulary \up\ is empty. 
 
It remains an open problem whether there is any logic at all (in the 
sense of \cite{logic}) capturing PTime on unordered structures.  It 
was conjectured in \cite{logic} that there is no such logic. 
 
There have been, however, continuing efforts to find logical systems 
capturing at least large parts of PTime, if not all of it.  These 
efforts have looked primarily in two directions.%
\footnote{A third direction, studied by Gire and Hoang \cite{gh},
  involves a form of restricted nondeterminism.  This direction looks
  promising, but we do not address it in this paper.}
One direction 
involves adding to FP additional constructs, usually in the form of 
quantifiers, to permit the direct expression of certain easily 
computable properties of unordered structures, for example the 
property ``the universe has an even number of elements'' mentioned 
above.  The most popular of these extensions has been to add counting 
to the logic.  There are several ways to formalize this extension; we 
choose the one described in \cite[Ch.~4]{otto}.  It involves adjoining 
to the input structure a second sort, consisting of the natural 
numbers up to the cardinality of the input set, and adding to the 
language terms of the form ``the number of elements $x$ satisfying 
$\phi(x)$.''    
 
The second direction taken by the search for a PTime logic involves 
combining a standard computation mechanism with additional logical
(rather 
than arithmetical) facilities.  The relational machines of Abiteboul
and 
Vianu \cite{av} are of this sort, combining a Turing machine and 
first-order logic.  Another model of this sort, more directly relevant
to 
our purposes here, is choiceless polynomial time, \cpt, introduced in 
\cite{bgs}.  Here the abstract state machine model of computation 
\cite{lipari, asm} is applied in a set-theoretic context, allowing 
essentially arbitrary data types over the input structure.  It is shown
in 
\cite{bgv} that \cpt\ is strictly stronger than PTime relational
machines, 
but even so it cannot compute the parity of an unstructured set 
\cite{bgs}.  It thus appears that this second direction produces unduly 
weak models.  On the other hand, we shall see that \cpt\ is capable of 
computing some things that are beyond the reach of FP plus counting. 
 
It is therefore reasonable to combine the two directions and consider 
computation models (or logics) like \cptc, which is \cpt\ augmented 
with the ability to compute cardinalities.  This model was already 
proposed in \cite[Subsection~4.8]{bgs} as worthy of further study. 
The present paper contains the first results of that study.  The main 
problem, which remains open, is whether \cpt\ plus counting captures 
polynomial time on unordered structures.    
 
Most of the results we present here are concerned with specific
algorithmic problems that are solvable in PTime but appear to be at
the borderline of expressibility in logics like \cpt\ plus counting.
Several of them are plausible candidates for separating PTime from
\cpt\ plus counting.
 
We begin with work motivated by the result from \cite{bgs} that 
bipartite matching is not in \cpt.  The proof of this involved 
exceptionally simple instances of the bipartite matching problem.  In 
the traditional picture of bipartite matching, where the input 
consists of a set of boys, a set of girls, and a (symmetric) ``willing 
to marry'' relation between them, the instances used in \cite{bgs} can 
be described as follows.  First suppose there are $2n$ boys and $2n$ 
girls, divided into two gangs of $n$ boys and $n$ girls each; a boy 
and a girl are willing to marry if and only if they belong to the same 
gang.  Obviously, a complete matching exists in this case.  Next 
suppose one of the boys defects from his gang and joins the other, 
while all girls remain in their original gangs.  Obviously there is no 
matching now.  But a \cpt\ program cannot distinguish these two 
situations if $n$ is sufficiently large compared to the program 
(\cite[Thm.~43]{bgs}).  This specific deficiency can evidently be 
removed by adding to \cpt\ the ability to count, but it seems that 
this success depends on the very simple structure of the ``willing to 
marry'' relation.  There seems to be no way to extend this result to 
general instances of bipartite matching.  Thus, there was some hope 
that bipartite matching would serve to separate PTime from \cptc. 
That hope is dashed here in Section~\ref{bipartite}, where we present 
a \cptc\ algorithm to decide whether a bipartite graph has a complete 
matching.  In fact, we show the rather surprising result that the 
existence of a complete matching can be expressed in FP+Card. 
 
In an effort to separate PTime from \cptc, we next turn to the two
known types of examples separating PTime from FP+Card.  These examples
involve certain graphs defined by Cai, F\"urer, and Immerman
\cite{cfi} and structures called multipedes introduced by two of the
present authors \cite{multipede}.  For the reader's convenience, we
recapitulate the relevant information from \cite{cfi} and
\cite{multipede}.  Then we discuss how the constructions from these
papers lead naturally to queries that are in PTime but not in FP+Card.
We show that, for suitably padded versions of the examples from
\cite{cfi} and also for the so-called 4-multipedes of \cite{multipede}
(without padding), these queries are in \cpt\ even without counting.
Thus, these examples show that FP+Card does not include \cpt\ and is
strictly included in \cptc.
 
There are very similar queries, using the graphs from \cite{cfi} 
without padding or using the 3-multipedes from \cite{multipede}, which 
are still in PTime (by somewhat trickier proofs than the versions in 
the preceding paragraph) but which we do not see how to express in 
\cptc.  So perhaps one of these will give the desired separation. 
 
Finally, motivated by the use of linear algebra modulo 2 in some of 
the preceding arguments, we consider the computation of (suitably 
presented) determinants.  We show that the question whether a matrix 
over a finite field or over the integers is singular (i.e., has zero 
determinant) is in FP+Card.  It is not in \cpt, even over the 
two-element field; the proof of this uses the zero-one law proved by 
one of us in \cite{sh634} and discussed by the other two in 
\cite{01law}.  The computation of determinants (in contrast to merely 
deciding whether they are zero) over the prime field $\bbb Z/p$ for an 
odd prime $p$ is in PTime, but we do not know whether it is in \cptc. 

\section{Background}		\label{background} 
 
In this section we review the logics and the computation models 
relevant to this paper, namely 
\begin{ls} 
\item fixpoint logic (FP), 
\item finite variable infinitary logic ($L^\omega_{\infty,\omega}$), 
\item choiceless polynomial time (\cpt), 
\end{ls}   
and their extensions by counting, FP+Card, $C^\omega_{\infty,\omega}$, 
and \cptc, respectively.  (The notation ``+Card'' stands for adding 
cardinality to the logic.)  We refer the reader to \cite{otto} for 
details about FP, FP+Card, $L^\omega_{\infty,\omega}$, and 
$C^\omega_{\infty,\omega}$ and to \cite{bgs} for details about \cpt\ 
and \cptc. 
 
Fixpoint logic FP is obtained by adding to ordinary first-order logic 
the (inflationary) fixpoint operator defined as follows.  If $X$ is an 
$r$-ary relation symbol not in \up, if $\phi(X,\vec x)$ is a formula 
of the vocabulary $\up\cup\{X\}$, and if $\vec x$ is an $r$-tuple of 
distinct variables, then $\text{FP}_{X,\vec x}\phi$ is used as an 
$r$-ary relation symbol.  It is interpreted as the fixed point 
obtained by starting with the empty relation and iterating the 
operation   
$$ 
R\mapsto R\cup\{\vec x:\phi(R,\vec x)\}. 
$$ 
(For a precise formulation, one should fix an \up-structure and values 
for all free variables of $\phi$ except $\vec x$, and then the 
operation above should be interpreted using these data; see 
\cite{otto}.)      
 
Fixpoint logic with counting, called FP+C in \cite{otto} but FP+Card 
here (to conform with the notation \cptc), is obtained from FP by the 
following modifications.  First, every input structure \ger A, with 
underlying set $A$, is replaced by a two-sorted structure $\ger A^*$ 
in which one sort is \ger A and the other is the initial segment 
$\{0,1,\dots,|A|\}$ of the natural numbers with the successor 
function.  Second, for each variable $x$ and formula $\phi(x)$, there 
is a term $(\sharp x)\phi(x)$ denoting the number (an element of the 
new sort) of values of $x$ that satisfy $\phi(x)$.  Fixpoint 
operations are allowed to define relations on either or both sorts; in 
particular, addition and multiplication are definable on the number 
sort insofar as their values don't overflow the available range of 
numbers.    
 
The infinitary logic $L^\omega_{\infty,\omega}$ is obtained from 
ordinary first-order logic by making two changes.  First, allow 
conjunctions and disjunctions of arbitrary, possibly infinite sets of 
formulas.  (The logic resulting from this first change is called 
$L_{\infty,\omega}$.)  Second, require each formula to use only a 
finite number of variables, where both free and bound variables are 
counted but the same variable may be re-used, i.e., it may occur both 
free and bound and possibly bound many times.  $L^k_{\infty,\omega}$ 
is the sub-logic in which the number of variables in a formula is 
required to be at most the natural number $k$.  It is known (see for 
example \cite[Cor.~1.30]{otto}) that FP is included in 
$L^\omega_{\infty,\omega}$ in the sense that, for every formula of FP, 
there is a formula of $L^\omega_{\infty,\omega}$ that is semantically 
equivalent, i.e., satisfied by the same tuples of elements in the same 
structures. 
 
The logic $L^\omega_{\infty,\omega}$ could be extended by counting 
terms just as FP was extended to FP+Card, but we shall instead follow 
\cite{otto} and use the more traditional counting quantifiers.  The 
logic $C^\omega_{\infty,\omega}$ is obtained from 
$L^\omega_{\infty,\omega}$ by adding the quantifiers $\exists^{\geq 
m}$ for all natural numbers $m$, semantically interpreted as ``for at 
least $m$ values of''.  It is shown in \cite[Cor.~4.20]{otto} that 
$C^\omega_{\infty,\omega}$ includes FP+Card. 
 
Since \cpt\ is newer and less widely known than the fixpoint and 
infinitary logics discussed above, we describe it in somewhat more 
detail, but for a full definition we refer to \cite{bgs}.  \cpt\ is 
the polynomial time fragment of a programming language, BGS, defined 
as follows.  Inputs to a computation are finite structures for a 
vocabulary \up; each program is associated with a fixed \up, but 
different programs can use different \up's and thus admit different 
sorts of inputs.  A computation proceeds in discrete stages, the state 
at any moment being a structure of the following sort.  Its underlying 
set \HF I consists of the underlying set $I$ of the input structure 
(regarded as a set of atoms, i.e., non-sets) plus all hereditarily 
finite sets over $I$, that is, all subsets of $I$, all sets whose 
members are either such subsets or members of $I$, etc.  In other 
words, \HF I is the smallest set having among its members all its 
finite subsets and all the members of $I$.  Notice that \HF I contains 
the natural numbers, coded as von~Neumann ordinals, 
$$ 
0=\emp, \quad 
1=\{\emp\}, \quad 
\dots,\quad 
n=\{0,1,\dots,n-1\},\quad\dots. 
$$ 
For computational purposes, this representation of the natural numbers
is
equivalent (in the BGS context) to unary notation for the natural
numbers. So we assume from now on that natural numbers (and in fact all
integers) are available, in unary notation, along with the basic
arithmetical operations.  In Section~\ref{det}, we shall also need
binary
notations; details about that representation will be given there.

We use 0 and 1 to represent the truth values \ttt{false} and 
\ttt{true}, respectively.  Thus predicates can be regarded as 
$\{0,1\}$-valued functions.  The structure giving a state of the 
computation has the following basic functions: 
\begin{ls} 
\item the functions and relations of the input structure, relations 
being regarded as $\{0,1\}$-valued functions, and all functions being 
extended to have value 0 when any input is not in $I$, 
\item the logical functions: $=,\ttt{true}, \ttt{false}, \neg, 
\land,\lor$, 
\item the set-theoretic functions $\in, \emp, \text{Atoms}, \bigcup, 
\text{TheUnique}, \text{Pair}$, 
\item finitely many dynamic functions, including \ttt{Halt} and 
\ttt{Output}. 
\end{ls} 
Here Atoms means the set $I$ of atoms (as opposed to sets) in \HF I. 
The function $\bigcup$ sends any $x$ to the union of all the sets that 
are members of $x$, $\text{TheUnique}(x)$ is the unique member of $x$ 
if $x$ is a set having exactly one member (and 0 otherwise), and 
$\text{Pair}(x,y)$ is the set $\{x,y\}$.  The dynamic functions are 
constant with value 0 in the computation's initial state but acquire 
more interesting values as the computation proceeds.  The vocabulary 
of a BGS program has symbols for all these functions.  The symbols for 
the input relations (as opposed to functions), the logic symbols, 
$\in$, the dynamic functions \ttt{Halt} and \ttt{Output}, and possibly 
some other dynamic function symbols are called \emph{Boolean}; their 
values are always 0 or 1.  (If we were interested in computing results 
other than ``yes'' and ``no'', then we would not declare \ttt{Output} 
to be Boolean.) 
 
The meaningful expressions of the programming language BGS are terms 
and rules.  Terms are built from the function symbols described above 
and variables in the usual way, with the addition of the term-forming 
construction   
$$ 
\{t(v):v\in r:\phi(v)\}, 
$$ 
where $t$ and $r$ are terms, $\phi$ is a Boolean term (i.e., one whose 
outermost constructor is a Boolean function symbol), and $v$ is a 
variable not free in $r$.  (By writing $v$ in the contexts $t(v)$ and 
$\phi(v)$, we mean to indicate only that $v$ is allowed to occur free 
there, not that it must occur free, nor that other variables cannot 
occur free.)  The interpretation of the term $\{t(v):v\in r:\phi(v)\}$ 
is the set of values of $t(v)$ for all values of $v$ that are members 
of the value of $r$ and make $\phi$ true.  When $\phi$ is \ttt{true},
we sometimes omit it from the notation and write simply $\{t(v):v\in
r\}$. 

We note for future use that there are terms representing the union of
two sets,
$$
a\cup b = \bigcup\text{Pair}(a,b)
$$
and the traditional set-theoretic coding of ordered pairs
$$
\sq{a,b}=\{\{a\},\{a,b\}\}=
\text{Pair}(\text{Pair}(a,a),\text{Pair}(a,b)).
$$
 
Rules are built by the following inductive construction.  Each rule 
defines, in an obvious way, a set of \emph{updates} of the state, 
provided values are specified for the rule's free variables. 
\begin{ls} 
\item \ttt{Skip} is a rule (producing no updates). 
\item If $f$ is a dynamic function symbol, say $j$-ary, and 
$t_0,t_1,\dots,t_j$ are terms, with $t_0$ Boolean if $f$ is, then   
$$ 
f(t_1,\dots,t_j):=t_0 
$$ 
is a rule. 
\item If $R_0$ and $R_1$ are rules and $\phi$ is a Boolean term, then   
$$ 
\ttt{if\ }\phi\ttt{\ then\ }R_0\ttt{\ else\ }R_1\ttt{\ endif} 
$$ 
is a rule. 
\item If $R_0(v)$ is a rule, $v$ is a variable, and $r$ is a term in 
which $v$ is not free, then   
$$ 
\ttt{do\ forall\ }v\in r, R_0(v)\ttt{\ enddo} 
$$ 
is a rule. 
\end{ls} 
The notion of free variable, used in these definitions, is defined in 
the usual way, with the term constructor $\{t(v):v\in r:\phi(v)\}$ and 
the rule constructor $\ttt{do\ forall\ }v\in r, R_0(v)\ttt{\ enddo}$ 
binding the variable $v$. 
 
A \emph{program} is a rule with no free variables.  To \emph{fire} a 
program in a state is to modify the dynamic functions of the state 
according to all the updates produced by the program \emph{except} 
that, if two of these updates are contradictory (i.e., update the same 
dynamic function at the same tuple of arguments to different values), 
then none of the updates are executed.  A \emph{run} of a program on 
an input is a sequence of states in which the first state is the 
initial state determined by the input structure (as above, with all 
dynamic functions constantly 0) and each subsequent state is obtained 
from its predecessor by firing the program.  The result of the 
computation is the value of \ttt{Output} at the first stage where 
\ttt{Halt} has the value \ttt{true} (i.e., 1).  (It would do no harm 
to automatically stop all runs whenever \ttt{Halt} has the value 
\ttt{true} or to insist that programs produce no updates in this 
situation.)  If \ttt{Halt} never becomes \ttt{true} then the 
computation fails to produce an output. 
 
A \emph{PTime bounded BGS program} is a BGS program $\Pi$ together
with two polynomials $p(n)$ and $q(n)$.  A run of $(\Pi,p(n),q(n))$ on
input $I$ is a run of $\Pi$ consisting of at most $p(|I|)$ stages and
having at most $q(|I|)$ active elements.  We do not reproduce here the
definition of ``active'' from \cite{bgs} but remark that, roughly
speaking, an element of \HF I is active if it is either involved in an
update during the run or a member of something involved in an update,
or a member of a member, etc.

For the purposes of this paper, we define
\cpt\ as the class of Boolean queries decidable by PTime
bounded BGS programs.  (A broader definition, using a three-valued
logic to accommodate computations where \ttt{Halt} never becomes true,
was used in \cite{bgs}, but we will not need to use it here.)

We observe that \cpt\ includes the expressive power of first order
logic.
The propositional connectives were included among the basic functions on
\HF I, and the quantifiers over the input structure can be simulated
because $(\exists v \in I)\,\phi(v)$ is equivalent to
$$ 
0\in\{0:v\in\text{Atoms}:\phi(v)\}. 
$$ 
Furthermore, \cpt\ includes the expressive power of fixpoint logic, 
for the iteration defining a fixpoint can be simulated by the 
iteration involved in the notion of run.  In fact, it was shown in 
\cite[Thm.~20]{bgs} that \cpt\ can simulate the PTime relational 
machines of Abiteboul and Vianu \cite{av}; it is known that these can 
compute all FP-definable queries. 
 
We shall need several times the observation that \cpt\ includes all
PTime 
(and in fact exponential time) computations on sufficiently small 
parts of the input structure.  Specifically, if the input structure 
has a definable subset $P$ with $|P|!\leq|I|$, then a BGS program can 
first produce, in a parallel computation, $|P|!$ subprocesses each of 
which knows a linear ordering of $P$.  Then each of these subprocesses 
can run a PTime algorithm on its ordered version of $P$.  If the 
PTime algorithm produces the same answer for all orderings, then these 
subprocesses will all give \ttt{Output} that value, so the overall 
algorithm produces this answer.  And the inequality $|P|!\leq|I|$ 
implies that this is a PTime bounded BGS program, so the result of the 
computation is in \cpt.  For the details of this argument, see 
the proof of \cite[Thm.~21]{bgs}.  Similarly, under the weaker 
assumption that $2^{|P|}\leq|I|$, a BGS program can produce, in 
polynomial time, all the subsets of $P$. 
 
To add counting to \cpt, we simply include, in every state, the 
additional function Card that sends every set to its cardinality 
(considered as a von~Neumann ordinal) and sends atoms to 0.  The 
resulting complexity class is called \cptc. 

The cardinality function makes it possible to carry out, in a single
step, the operations of addition and multiplication on von~Neumann
ordinals.  Indeed, we can express $a+b$ as the cardinality of
$$
a\cup\{\sq{0,x}:x\in b\},
$$
and $ab$ is the cardinality of the cartesian product
$$
a\times b = \bigcup\{\{\sq{x,y}:x\in a\}:y\in b\}.
$$

\begin{rmk}
Theorem~8 of \cite{bgs} says, roughly speaking, that every object
activated during a run of a PTime BGS program was ``looked at'' during
that run.  This is no longer true when we add Card to the computation
model.  The ordinal $\text{Card}(x)$ can be active in a state without
all its predecessors being looked at.  For example, if the input is a
linearly ordered set of size $n$, a computation can, since addition is
available, initialize a nullary dynamic $c$ to 1 and then perform $n$
steps doubling $c$ at each step.  Then $2^n$ is active in the final
state, but most of the ordinals below it have not been looked at. 

Intuitively, the computation just described should not count as
polynomial
time, for it involves parallel processes indexed by sets of size
exponentially big compared to the input.  Our definition of PTime in the
BGS context agrees with this intuition, for the number of active
elements
in this computation is exponential.  (The number of critical elements,
in
the sense of \cite{bgs} is only polynomial, so it is important to
include
members of critical elements, their members, etc.\ in the definition of
active elements and thus in the definition of PTime.
\end{rmk}

\begin{rmk}
BGS was designed for theoretical purposes.  Some of its conventions
were designed to simplify analysis of programs and thus are unnatural
from a programming point of view.  In this paper, we retain those
conventions and work around them where necessary.  But for actual
programming, these conventions should be modified.  In particular,
arithmetic should be available directly rather than being coded in
the von~Neumann ordinals.  The input itself should in general be a
metafinite structure in the sense of \cite{gg}.
\end{rmk}

The following diagram indicates the relationships between the various 
logics and complexity classes considered here.  Arrows represent 
inclusion relationships. 
$$ 
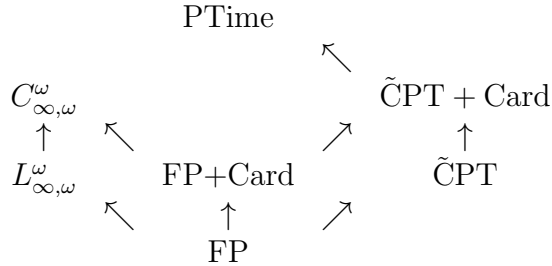
\begin{array}{ccccc} 
&&\text{PTime}&&\\ 
&&&\nwarrow&\\ 
C^\omega_{\infty,\omega}&&&&\cptc\\ 
\uparrow&\nwarrow&&\nearrow&\uparrow\\ 
L^\omega_{\infty,\omega}&&\text{FP+Card}&&\cpt\\ 
&\nwarrow&\uparrow&\nearrow&\\ 
&&\text{FP}&& 
\end{array} 
$$ 
In the following sections, we consider various specific problems and 
try to determine which of these logics and complexity classes contain 
them.   
 
\section{Bipartite Matching}	\label{bipartite} 
 
\subsection{Statement of the Problem} 	\label{match}   
   
\begin{df}   
A \emph{bipartite graph} consists of two finite sets $A$ and $B$ with   
an adjacency relation $R\subseteq A\times B$.   
\end{df}   
   
We denote a bipartite graph by \sq{A,B,R}.  It makes no difference   
whether we regard it as a two-sorted structure with one binary   
predicate $R$ or as a one-sorted structure with, in addition to $R$,   
unary predicates for $A$ and $B$.  Even if we adopt the two-sorted   
viewpoint, we assume whenever convenient that $A$ and $B$ are   
disjoint.   
   
\begin{df}   
A \emph{matching} in a bipartite graph \sq{A,B,R} is a partial   
one-to-one function $M$ from $A$ into $B$ which, considered as a   
binary relation (a set of ordered pairs), is a subset of $R$.  We call

a matching \emph{complete} if its   
domain is all of $A$.   
\end{df}   
 
\begin{df}   
\emph{Bipartite matching} is the following decision problem.  The   
input is a bipartite graph \sq{A,B,R} and the question is whether it   
has a complete matching.   
\end{df}   
 
In these definitions, we did not require $|A|=|B|$, so a complete 
matching may have its range strictly included in $B$.  Everything we 
say about the bipartite matching problem remains true if we restrict 
the inputs to be bipartite graphs \sq{A,B,R} with $|A|=|B|$. 
   
It was shown in \cite{bgs} that the bipartite matching problem is not 
in \cpt.  The proof of this fact exploited the inability of \cpt\ to 
count.  For the particular graphs used in that proof, the decision 
problem would become easy if counting were available, but this 
observation does not apply to more general instances of bipartite 
matching.  The question thus arises whether the bipartite matching 
problem is in \cptc. 
 
We present here a somewhat surprising affirmative answer.  In fact, we 
show that bipartite matching is expressible in FP+Card. 
   
\subsection{Known Algorithms} 	\label{known} 
   
In this subsection we describe two well-known approaches to the
bipartite matching problem.  Neither provides a solution in \cptc, but
both will play a role in the solution.
   
The first is an algorithm which we call the \emph{path algorithm}.  It
works with (incomplete) matchings, starting with the empty one and at   
each stage either replacing the current matching by a larger one or   
determining that no complete matching can exist.   
   
To describe a step of this algorithm, let $M$ be the current matching.
If it is complete, then output ``yes'' and halt.  If it is incomplete,
then proceed as follows.  Consider the directed graph whose vertex set
is $A\dot\cup B$ (we invoke our standing assumption that $A$ and $B$   
are disjoint whenever convenient) and whose directed edges are   
\begin{ls}   
\item all $(a,b)\in R-M$ and   
\item the converses $(b,a)$ of all $(a,b)\in M$.   
\end{ls}   
In other words, start by regarding all pairs $(a,b)\in R$ as directed   
edges from $a$ to $b$, but then reverse the direction of those pairs   
that are in the current matching $M$.     
   
If this directed graph has a directed path from a vertex $a\in   
A-\dom M$ to a vertex $b\in B-\ran M$, then choose one such path, and   
let $P$ be the corresponding set of pairs in $R$ (i.e., take the edges
in the path and reverse the direction of those that are in $M$, so as   
to get pairs in $R$).  Notice that, except for the endpoints $a$ and   
$b$, every vertex in our path has two incident edges in $P$, one of   
which is in $M$ and the other not in $M$; the endpoints, of course,    
have only one incident edge (each) and it is not in $M$.  This implies
immediately that the symmetric difference $M\bigtriangleup P$ is a   
matching of cardinality one greater than that of $M$.  Proceed to the   
next step with $M\bigtriangleup P$ as the current matching.   
   
If the directed graph has no path from a vertex $a\in A-\dom M$ to a   
vertex $b\in B-\ran M$, then output ``no'' and halt.   
   
This completes the description of the algorithm, but it should be   
accompanied by an explanation of why the ``no'' answer in the last   
situation is correct.  (All other aspects of correctness --- eventual   
termination and correctness of the ``yes'' answers --- are obvious.)   
So suppose that, at some step of the algorithm, there is no directed   
path of the required sort.  Since the algorithm has not yet halted   
with ``yes'', there are points $a\in A-\dom M$; fix one such $a$.  Let
$X$ and $Y$ be the sets of all vertices in $A$ and $B$, respectively,   
that are reachable from $a$ by directed paths in    
the digraph under  consideration.  By assumption, $Y\subseteq\ran M$.   
Furthermore, by the definition of the digraph, $X$ contains all the   
points that are matched by $M$ with points in $Y$.  In addition, $X$   
contains the point $a$, which isn't matched with anything.  Therefore,
$|X|>|Y|$.  It is easy to check, using again the definition of the   
digraph, that every pair in $R$ whose first component is in $X$ has   
its second component in $Y$.  Thus, if there were a complete matching   
for \sq{A,B,R}, it would have to map $X$ one-to-one into $Y$, which is
impossible as $|X|>|Y|$.  Therefore, no such matching can exist, and   
the algorithm's ``no'' answer is correct.   
   
The preceding discussion not only establishes the correctness of the 
path algorithm but also essentially proves the celebrated ``marriage 
theorem'', often called Hall's theorem.  For the history of this 
theorem see \cite[page 54]{bol}. 
   
\begin{thm}[``Marriage'']   
A bipartite graph \sq{A,B,R} admits a complete matching if and only if,
for every $X\subseteq A$, its set of $R$-neighbors    
$$   
Y=\{b\in B:(\exists x\in X)\,(x,b)\in R\}   
$$   
has cardinality $|Y|\geq|X|$.   
\end{thm}   
   
\begin{pf}   
As mentioned above, a complete matching must map every $X$   
one-to-one into the corresponding $Y$, so $|Y|\geq|X|$.   This proves
the ``only if'' part of the theorem.
 
For the ``if'' part, suppose \sq{A,B,R} has no complete matching.
Then the path algorithm must eventually output ``no,'' and when it
does it has found, according to the discussion above, an $X$ (namely
the set of points reachable in the digraph from an $a\in A-\dom M$)
for which the corresponding $Y$ has $|Y|<|X|$.
\end{pf}   
   
A second approach to the bipartite matching problem would be to check   
the condition in the marriage theorem.   
This second approach is clearly not in polynomial time, for there are   
exponentially many $X$'s to check.  It is choiceless, as the   
computations for all $X$'s can be done in parallel (after all the   
$X$'s have been generated, in another parallel computation), but, as   
pointed out in \cite{bgs}, choicelessness is hardly relevant in the   
absence of a bound on the number of activated sets.  (There is a   
little bit of relevance, since to simulate choice by parallel   
computation one generally needs to activate $n!$ sets, where $n$ is   
the input size, and the algorithm based on the marriage theorem   
activates only approximately $2^n$ sets.)   
   
The path algorithm, by contrast, clearly runs in polynomial time, 
since it will terminate after at most $n$ steps and each step consists 
mainly of testing the existence of paths between certain vertices, 
which can be done in polynomial time.  Unfortunately, this algorithm 
requires arbitrary choices.  Whether the required path exists at any 
step can be decided choicelessly (see \cite[Section~1]{bgs}), but the 
algorithm requires choosing one such path in order to form the 
matching $M\bigtriangleup P$ for the next step to use. 
   
In fact, since the path algorithm not only decides whether a complete   
matching exists but, when the answer is ``yes,'' produces one, it   
clearly cannot be choiceless, for some bipartite graphs have complete   
matchings but none that are invariant under all the graph's   
automorphisms.  For example, consider the case where $|A|=|B|\geq2$   
and $R=A\times B$.   
   
We record for future reference that, if the input is given with a   
linear order (i.e., if each of $A$ and $B$ is linearly ordered), then   
no further choices are needed by the path algorithm.  When it searches
for a path, it can do a depth-first search, going through vertices in   
the given order, and use the first path that it finds.   
   
\subsection{A Choiceless Polynomial Time Algorithm}	\label{alg}   
   
In this subsection, we describe a \cptc\ algorithm that solves the 
bipartite matching problem.  The description will be informal, but it 
should be clear that the algorithm could be programmed in BGS 
augmented by the cardinality function Card and that it would run in 
polynomial time.  Thus, it shows that bipartite matching is in \cptc. 
Later, we shall prove, somewhat more formally, that the algorithm 
enables us to express ``there is a complete matching'' in the language 
FP+Card; this implies formally that the problem is in \cptc. 
   
The algorithm proceeds in three phases, given a bipartite graph
\sq{A,B,R} as input.     
   
In phase 1, we partition $A$ and $B$ into subsets $A_i$ ($i\in I$) and
$B_j$ ($j\in J$) respectively in such a way that    
\begin{ls}   
\item for each $i$ and each $j$, all the vertices in $A_i$ have the   
same number of $R$-edges to $B_j$.  That is,    
$$   
(\forall a,a'\in A_i)\quad|\{b\in B_j:(a,b)\in R\}|=    
|\{b\in B_j:(a',b)\in R\}|,   
$$   
\item symmetrically,   
$$   
(\forall b,b'\in B_j)\quad|\{a\in A_i:(a,b)\in R\}|=    
|\{a\in A_i:(a,b')\in R\}|,    
$$   
and   
\item the index sets $I$ and $J$ have canonical linear orderings.   
\end{ls}   
This is achieved by the following procedure, called the ``stable   
coloring algorithm'' in \cite{otto}.   
   
We proceed in steps, having at each step a partition of $A$ and a   
partition of $B$, together with, for each of these partitions, a   
linear ordering of the blocks.  As long as the current partitions are   
not of the desired sort, they will be refined, i.e., replaced with new
partitions each block of which is included in a block of the   
corresponding old partition.  We begin with each of $A$ and $B$   
trivially partitioned into a single piece (so there is no question   
about the linear ordering of blocks).   
   
Suppose, at some stage, we have a partition that does not satisfy the   
requirements listed above.  Since we have linear orderings as in the   
third requirement, there must be a violation of one or both of the   
first two requirements.     
   
Replace each block $A_i$ by a sequence of subblocks determined as   
follows.  To each $a\in A_i$ assign a vector consisting of the   
cardinalities $|\{b\in B_j:(a,b)\in R\}|$ listed in order of the   
blocks $B_j$.  Each vector so obtained will give one subblock,   
consisting of all the $a\in A_i$ that produced that vector.  The   
subblocks within $A_i$ are ordered according to the lexicographic   
ordering of their vectors.  Subblocks coming from different blocks   
$A_i$ and $A_{i'}$ are ordered as those blocks were ordered in the   
given partition.   
   
Replace each block $B_j$ by subblocks and linearly order these 
subblocks analogously. 
   
Since at least one of the first two requirements was violated, at   
least one of our two partitions will be properly refined.  Thus, the   
number of steps of this sort is bounded by (slightly less than) the   
number of vertices in $A\dot\cup B$.  Therefore, phase 1 must   
terminate, and when it does it has provided partitions satisfying all   
our requirements.   
   
In phase 2 we first replace $R$ by the following (possibly) larger   
relation:    
$$   
R^+=\bigcup\{A_i\times B_j:(A_i\times B_j)\cap R\neq\varnothing\}.   
$$   
In other words, as soon as one vertex in $A_i$ is $R$-joined to one   
vertex in $B_j$ (and therefore, by the requirements on the partition,   
every vertex in $A_i$ is joined to at least one vertex in $B_j$, and   
vice versa), $R^+$ joins every vertex in $A_i$ to every vertex in   
$B_j$.     
   
Then, using the linear ordering of the blocks produced in phase 1, we   
create an isomorphic copy \sq{A',B',R'} of \sq{A,B,R^+} in which the   
vertex sets $A'$ and $B'$ are equipped with canonical linear   
orderings.  To do this, let $A'$ consist of triples $(0,i,r)$ where
$i\in   
I$ and $r$ is a natural number in the range $0\leq r<|A_i|$.  The idea
is that the $|A_i|$ triples whose second component is $i$ act as a   
substitute for the members of $A_i$.  Define $B'$ analogously, using   
triples $(1,j,s)$, and let   
\begin{eqs}   
((0,i,r),(1,j,s))\in R'&\iff& (a,b)\in R^+\text{ for some (all) }a\in   
A_i,\,b\in B_j\\   
&\iff& (a,b)\in R\text{ for some }a\in A_i,\,b\in B_j.   
\end{eqs}%
Notice that $\sq{A',B',R'}$ is isomorphic to \sq{A,B,R^+}, but there   
is no canonical choice of an isomorphism.  To choose a specific   
isomorphism we would need to choose a linear ordering of each of the   
sets $A_i$ and $B_j$.  Fortunately, the algorithm doesn't need any   
isomorphism, so it remains choiceless.   
   
Finally, in phase 3, we apply the path algorithm to determine whether   
\sq{A',B',R'} has a complete matching.  No arbitrary choices are   
involved here, since $A'$ and $B'$ are (unlike $A$ and $B$)   
canonically linearly ordered: Their elements are triples whose second   
components come from the index sets $I$ and $J$, for which phase 1   
provided a linear order, and whose first and third components are   
natural numbers. So we can use the lexicographic order on the   
triples.    
   
Output ``yes'' or ``no'' according to whether \sq{A',B',R'} has a   
complete matching or not.   
   
This completes the description of the algorithm.  It should be clear   
that it is in \cptc, but there is a real question about its   
correctness.  The next subsection addresses that question.   
   
\subsection{Correctness Proof}	\label{correct}   
   
If the algorithm presented in the last subsection outputs ``no,'' this
means that there is no complete matching in \sq{A',B',R'}, hence no   
complete matching in the isomorphic graph \sq{A,B,R^+}, and hence no   
complete matching in the original graph \sq{A,B,R}, because   
$R\subseteq R^+$ so any complete matching for $R$ would also be one for
$R^+$.     
   
If, on the other hand, the algorithm ouputs ``yes,'' then there is a   
complete matching for \sq{A,B,R^+}, but this need not be a complete   
matching for \sq{A,B,R}, since it could use edges from $R^+-R$.  Thus,
the following lemma is needed to establish the correctness of the   
algorithm.     
   
\begin{la}   
In the situation of the preceding subsection, if \sq{A,B,R^+} has a   
complete matching, then so does \sq{A,B,R}.   
\end{la}   
   
\begin{pf}   
Fix a complete matching $M$ for \sq{A,B,R^+}.   
   
We first define a (reasonably fair) allocation of responsibility,    
among the edges $(a,b)\in R$, for the edges $(p,q)\in R^+$.  By   
definition of $R^+$, the fact that it contains $(p,q)$ is caused by   
the presence in $R$ of some edges $(a,b)$ between the same blocks   
$A_i$ and $B_j$.  If the number of such edges in $R$ is $n_{ij}$, then
we allocate responsibility for $(p,q)$ equally among them, assigning   
each such $(a,b)$ the amount $1/n_{ij}$ of responsibility.  Thus, the   
total responsibility (of all $(a,b)$) for one edge $(p,q)$ is 1, and   
the responsibility is shared by the $R$-edges between the same blocks   
as $(p,q)$.   
   
If $(p,q)$ is an edge in $R^+$ between blocks $A_i$ and $B_j$ and
$(a,b)$   
is an edge in $R$ but not between blocks $A_i$ and $B_j$ then the   
responsibility of $(a,b)$ for $(p,q)$ is zero.  If $S$ is a subset of
$R$   
then the responsibility of $S$ for an edge $(p,q)\in R^+$ is the sum of
the responsibilities of the edges in $S$ for $(p,q)$.  The
responsibility   
of $S$ for a subset of $R^+$ is the sum of the responsibilities of $S$
for   
the edges in the subset.  Further, a vertex $v$ gives rise to a subset
$S(v)$ of $R$, namely the set of edges of $R$ incident to that vertex.
A   
set $V$ of vertices also gives rise to a subset of $R$, namely the union
of the sets $S(v)$ where $v$ ranges over $V$.  The allows us to speak   
about the responsibility of a vertex or a set of vertices for an edge or
a   
set of edges in $R^+$.   
   
Because all vertices in $A_i$ have the same number of $R$-edges to
$B_j$,   
they all have equal responsibility for any $(p,q)$ joining these blocks
in   
$R^+$, namely responsibility $1/|A_i|$.  (In more detail: Each of these
vertices is incident to $n_{ij}/|A_i|$ edges to $B_j$, and each of these
edges bears responsibility $1/n_{ij}$ for $(p,q)$.  So each vertex has
responsibility $1/|A_i|$ for $(p,q)$.)   
   
Let $X$ be an arbitrary subset of $A$, and let $Y$ be, as in the
marriage   
theorem, the set of all vertices in $B$ that have an edge in $R$ from
some   
vertex in $X$.  We shall prove that $|X|\leq|Y|$; then the marriage   
theorem will provide the required matching for \sq{A,B,R}.   
   
Temporarily restrict attention to one block $A_i$.  We consider the   
total responsibility of vertices in $X\cap A_i$ for the edges of $M$   
(the fixed matching for $R^+$) that connect $A_i$ with $B$.  As noted   
above, each vertex in $A_i$ has the same responsibility $1/|A_i|$ for   
each such edge, so the vertices in $X\cap A_i$ have proportionate   
responsibility $|X\cap A_i|/|A_i|$ for each such edge of $M$.  There   
are, since $M$ is a complete matching, exactly $|A_i|$ such edges.   
Therefore, the total responsibility of all vertices in $X\cap A_i$ for
edges in $M$ (from $A_i$ to $B$) is $|X\cap A_i|$.  We wrote ``from   
$A_i$ to $B$'' in parentheses, because it can safely be omitted;   
vertices in $A_i$ have, by definition, no responsibility for edges (of
$M$ or otherwise) originating in other blocks $A_{i'}$.     
   
Now consider all the blocks, and sum over $i$ the result of the   
preceding paragraph.  The total responsibility of all vertices in $X$   
for all edges in $M$ is exactly $|X|$.   
Recalling how responsibility of vertices was defined, we can restate   
the result as: The total responsibility, for edges in $M$, of all   
$R$-edges originating in $X$ is exactly $|X|$.   
   
Now we repeat, as far as possible, the preceding three paragraphs   
``from the other side,'' i.e., starting with a fixed block $B_j$ in   
$B$ and computing the total responsibility of the vertices in $Y\cap   
B_j$ for the edges in $M$ (from $A$ to $B_j$).  It is their   
proportionate share, $|Y\cap B_j|/|B_j|$, of the total responsibility   
of $B_j$ for the $M$-edges that end in $B_j$.  So far, this is exactly
analogous to the preceding argument, but the next step is slightly   
different.  Although the domain of $M$ is, by completeness, all of   
$A$, its range need not be all of $B$.  Thus, the number of $M$-edges   
ending in $B_j$ is $\leq |B_j|$ (as $M$ is one-to-one), but equality   
need not hold.  Therefore, we can conclude only that the total   
responsibility of $Y\cap B_j$ for $M$-edges is $\leq|Y\cap B_j|$.     
   
Summing over all blocks $B_j$, we find that the total responsibility 
for $M$-edges of all vertices in $Y$ is $\leq |Y|$. 
As before, we rephrase this in terms of responsibility of edges: The   
total responsibility, for edges in $M$, of all $R$-edges ending   
in $Y$ is $\leq|Y|$.   
   
Finally, we recall that, by definition of $Y$, every $R$-edge   
originating in $X$ must end in $Y$.  Therefore   
\begin{eqs}   
|X|&=&\text{total responsibility for }M\text{-edges of }\\   
&&\text{edges originating in }X\\   
&\leq&\text{total responsibility for }M\text{-edges of }\\   
&&\text{edges ending in }Y\\   
&\leq&|Y|.   
\end{eqs}%
This completes the verification that \sq{A,B,R} satisfies the   
condition in the marriage theorem and therefore has a perfect   
matching.    
\end{pf}   
   
\subsection{Fixed-Point Logic With Counting} \label{fp+c}   
   
In this subsection we indicate how to express the existence of a 
complete matching in the extension FP+Card of first-order logic by 
fixed-point operators and counting. 
   
First, observe that what Otto calls the stable coloring in   
\cite[Section~2.2]{otto} amounts to the partitions $\{A_i:i\in I\}$   
and $\{B_j:j\in J\}$ produced by our algorithm together with the   
linear ordering of index sets construed as a pre-ordering of $A$ and   
$B$.  The pre-ordering has $x\prec y$ if and only if $x$ is in an   
earlier block than $y$.  (Technically, Otto works with one-sorted   
structures, so his stable coloring also has a convention for the   
relative ordering of $A$ and $B$, say $a\prec b$ whenever $a\in A$ and
$b\in B$.  This technicality will not matter in the following.)   
   
By \cite[Theorems~2.23 and 2.25]{otto}, the stable coloring is   
definable in the logic $C^2_{\infty,\omega}$, and its equivalence   
classes, our $A_i$'s and $B_j$'s, are exactly the equivalence classes   
with respect to $C^2_{\infty,\omega}$-equivalence.  Therefore, the
invariant   
$I_{C^2}$ as defined in \cite[Section~3.2]{otto} encodes all the
following    
information (plus more information that we won't need):   
\begin{ls}   
\item The blocks $A_i$ and $B_j$, regarded as points.   
\item The linear ordering of these blocks.   
\item For each pair of blocks $A_i$ and $B_j$ whether there is an   
$R$-edge joining them.     
\item The cardinality of each block.   
\end{ls}   
   
Recall from Section~\ref{background} that, for the logic FP+Card, 
structures \ger A (like our graphs \sq{A,B,R}) are enriched with a new 
sort containing the natural numbers from 0 up to and including the 
size of \ger A, with the standard successor function, and the 
resulting structure is called $\ger A^*$.  The standard linear 
ordering of natural numbers is easily FP-definable in $\ger A^*$.  
 
According to \cite[Lemma~4.14(ii)]{otto}, the invariant $I_{C^2}$ of
\ger   
A is FP+Card interpretable in $\ger A^*$, as a structure on the   
new, numerical sort.  This means that the linear ordering of the blocks
of   
$I_{C^2}$ is used to replace these blocks by numbers, the rest of the   
structure of $I_{C^2}$ is transferred to this copy, and the result is   
FP+Card-definable in $\ger A^*$.   
   
Once we have this form of $I_{C^2}$, we essentially have the structure
that we called \sq{A',B',R'} in our description of the algorithm in   
Subsection~\ref{alg}.  We can take $A'$ to be the set of pairs $(0,i,r)$
where $i$ is the number representing a block $A_i$ and $r<|A_i|$, and we
can take $B'$ to be the set of pairs $(1,j,s)$ where $j$ represents a   
block $B_j$ and $s<|B_j|$.  $R'$ joins $(0,i,r)$ to $(1,j,s)$ if there
was an $R$-edge from $A_i$ to $B_j$ --- information that we saw is  
available in $I_{C^2}$.     
   
Thus, we have a copy of \sq{A',B',R'}, with linearly ordered underlying
set (since it's in the numerical sort), FP+Card-definable in
$\sq{A,B,R}^*$.   
To apply the path algorithm to this copy is to apply a polynomial time
algorithm to an ordered structure.  So the result, the decision whether
there is a matching, is expressible in FP+Card, in fact in just FP,   
over this copy.  Therefore, the decision is expressible in FP+Card   
over $\sq{A,B,R}^*$, as claimed.   
 
\begin{rmk}
It is natural to ask whether, when a bipartite graph does not admit a
complete matching, one can compute in \cptc\ the size of the largest
(incomplete) matching.  

Our definition above (and in \cite{bgs}) of BGS programs allowed only
Boolean output, so technically one cannot compute in \cptc\ anything
other than Boolean queries.  But this restriction in the definition
was only a matter of convenience.  The BGS computation model and thus
the complexity classes \cpt\ and \cptc\ can and should be extended to
allow non-Boolean output whenever this is useful.

Once this extension is made, it is easy to show that the size of the
largest matching in a bipartite graph is computable in \cptc.  Indeed,
$\sq{A,B,R}$ has a matching whose domain contains all but $s$ elements
of $A$ if and only if there is a complete matching in the graph
obtained by adding $s$ new elements to $B$ and enlarging $R$ so as to
relate all elements of $A$ to the $s$ newly added elements.
\end{rmk}

The algorithm presented in this section depends on the fact that we deal
with a bipartite graph.  The notion of matching makes sense more
generally.  In any undirected, loopless graph, a \emph{matching} is a
family of edges no two of which have a common endpoint.  A matching is
\emph{complete} is every vertex is incident to an edge from the
matching.
Whether a given graph has a complete matching can be decided in
polynomial
time by a variant of the path algorithm.  But we do not know whether
this
decision can be computed without using choices or an ordering.

\begin{ques}
Is the existence of complete matchings in general (non-bipartite)
graphs computable in \cptc?
\end{ques}

\section{Cai-F\"urer-Immerman Graphs}	\label{cfi-graphs} 

In describing the Cai, F\"urer, Immerman construction, we follow, with 
a minor modification, Otto's presentation \cite[Example~2.7]{otto}, 
which is itself a minor modification of the presentation in 
\cite{cfi}. 
 
Let $G$ be a finite connected graph; we shall need only the special 
case where $G$ is the complete graph on some number $m+1$ of points, 
but the construction is no harder to present in the general case.  For 
each vertex $v$ of $G$, let $E_v$ be the set of edges incident with 
$v$.  Fix a linear ordering $\preceq$ of the vertices of $G$.  Using 
$G$, we define a new graph $G^*$ as follows.  Each vertex of degree 
$d$ in $G$ gives rise to $2^d$ vertices of $G^*$ and each edge of $G$ 
gives rise to two vertices of $G^*$.  Specifically, we let the 
vertices of $G^*$ be 
\begin{ls} 
\item pairs $(v,X)$ where $v$ is a vertex of $G$ and $X$ is a subset 
of $E_v$, and 
\item pairs $(e,+)$ and $(e,-)$ where $e$ is an edge of $G$.   
\end{ls} 
For each vertex $v$ of $G$, we write $U(v)$ for the set of associated 
vertices $(v,X)$ of $G^*$, and similarly for each edge $e$ of $G$, we 
write $U(e)$ for the pair of associated vertices $(e,\pm)$.  (For
vertices $(v,X)\in U(v)$, we chose to use, as members of the second
component $X$, edges $\{v,w\}\in E_v$ rather than simply the
distant vertices $w$ of those edges.  The main reason for this choice
is to match a visualization in which $U(v)$ is given in terms of
``local data'' at $v$ if edges are viewed as line segments.)

The edges of $G^*$ are also of two sorts; whenever edge $e$ and vertex
$v$ are incident in $G$, we
\begin{ls} 
\item join $(v,X)$ to $(e,+)$ if $e\in X$, and 
\item join $(v,X)$ to $(e,-)$ if $e\notin X$. 
\end{ls} 
In addition, we transport the linear ordering $\preceq$ of the 
vertices of $G$ to a pre-ordering, also called $\preceq$, on the 
vertices of the form $(v,X)$ in $G^*$; that is, we put 
$(v,X)\preceq(v',X')$ just in case $v\preceq v'$. 
 
Before proceeding with the construction, it is useful to analyze the 
automorphisms of the structure (graph with a linear pre-ordering of 
some of the vertices) $\ger G^*=\sq{G^*,\preceq}$.  Preserving 
$\preceq$, such an automorphism $\alpha$ must map each $U(v)$ into 
itself.  Then, to preserve adjacency, it must map each $U(e)$ into 
itself, for vertices in different $U(e)$'s have neighbors in different 
$U(v)$'s.  Thus, $\alpha$ gives rise to a subset $S$ of the edge set 
of $G$, namely 
$$ 
S=\{e:\alpha\text{ interchanges }(e,+)\text{ and }(e,-)\}. 
$$ 
Obviously, $S$ determines the action of $\alpha$ on vertices of the 
form $(e,\pm)$.  In fact, it completely determines $\alpha$, for a 
vertex $(v,X)$ is determined by which $(e,\pm)$'s are adjacent to it. 
More formally, we have 
$$ 
\alpha(e,\pm)= 
\begin{cases} 
(e,\mp)&\text{if }e\in S\\ 
(e,\pm)&\text{if }e\notin S 
\end{cases} 
\quad\text{and}\quad 
\alpha(v,X)=(v,X\bigtriangleup(S\cap E_v)). 
$$
(As before, $\bigtriangleup$ denotes symmetric difference.)  Conversely,
for any set $S$ of edges of $G$, the preceding formulas define an
automorphism of $\ger G^*$.
 
The Cai, F\"urer, Immerman graphs are subgraphs of $G^*$ obtained as 
follows.  For any subset $T$ of the vertex set of $G$, let $G^T$ be 
the induced subgraph of $G^*$ containing all the vertices of the form 
$(e,\pm)$ but containing $(v,X)$ only if either $v\in T$ and $|X|$ is 
odd or $v\notin T$ and $|X|$ is even.  An analysis exactly like that 
in the preceding paragraph shows that any automorphism of any $\ger 
G^T$ (meaning of course $\sq{G^T,\preceq}$) and in fact any 
isomorphism from one $\ger G^T$ to another $\ger G^{T'}$ must be given 
by the formulas above, for some set $S$ of edges of $G$.  To describe 
which $S$'s give isomorphisms between which $\ger G^T$'s, it is 
convenient to use the notation 
$$ 
\od S=\{v:\text{The number of edges in }S\text{ incident to }v 
\text{ is odd.}\}. 
$$ 
Then the $\alpha$ associated to $S$ maps $G^T$ to $G^{T'}$ just in 
case $T\bigtriangleup T'=\od S$ (equivalently, $T'=T\bigtriangleup\od 
S$). 
 
At this point we must recall two well-known facts from graph theory.
The first is that \od S always has even cardinality.  Indeed, the
total number of incident point-edge pairs $(v,e)$ is even because
every edge $e$ contributes two such pairs.  But, classifying the same
pairs according to their vertex components $v$, we find that the
number of these pairs is $\sum_v\text{degree}(v)$.  Modulo two, this
sum is congruent to the number of odd summands, i.e., to the
cardinality of \od S.  So this cardinality must, like the sum, be
even. 

The second fact to recall is that for connected graphs, like our $G$,
there is a converse to the first fact: Any set consisting of an even
number of vertices is \od S for some set $S$ of edges.  To see this,
pair off the vertices in the given set (arbitrarily) and choose for
each pair a path joining them.  Of course, if $P$ is the set of edges
of one of these paths, then \od P consists just of the two endpoints
of that path.  Summing up this information over all the chosen paths
$P$ and reducing modulo two, we find that our given set of vertices is
\od S where $S$ consists of those edges that occur in an odd number of
the paths $P$.

Applying these facts to our situation, we see that $\ger G^T$ and
$\ger G^{T'}$ are isomorphic if and only if $|T|$ and $|T'|$ have the
same parity.  We write $\ger G^0$ (resp.\ $\ger G^1$) for $\ger G^T$
when $|T|$ is even (resp.\ odd).  Thus, $\ger G^0$ and $\ger G^1$ are
well-defined up to isomorphism and are not isomorphic to each other.
 
In fact, $\ger G^0$ and $\ger G^1$ can be distinguished by the 
following simple property.  In $G^0$ it is possible to choose one from 
each pair of vertices $(e,\pm)$ corresponding to an edge of $G$, in 
such a way that each block $U(v)$ contains a vertex $(v,X)$ adjacent 
to the chosen vertices in all the pairs $U(e)$ for $e$ adjacent to $v$ 
in $G$.  This is easiest to see if we think of $G^0$ as $G^\emp$; then 
we simply choose $(e,-)$ from every $U(e)$.  (If we think of $G^0$ as 
$G^T$ for some other $T$ of even size, then we should fix an $S$ with 
$\od S=T$, and we should choose $(e,+)$ if and only if $e\in S$.)  On 
the other hand, no such choice is possible in $G^1$.  Indeed, let us 
represent $G^1$ as $G^T$ for a specific $T$ of odd size, and suppose a 
successful choice of $(e,\pm)$'s had been made.  Let $S$ be the set of 
edges $e$ where the choice was $(e,+)$.  Then we would have, for each 
vertex $v$ of $G$, that $(v, S\cap E_v)\in G^T$, which means that   
$$ 
v\in T\iff |S\cap E_v|\text{ is odd }\iff v\in\od S. 
$$ 
So $T=\od S$, which is absurd as $|T|$ is odd and $|\od S|$ is even. 
 
Let us now specialize to the case where $G$ is the complete graph on 
$m+1$ vertices.  Let $\ger H_m^0$ and $\ger H_m^1$ be padded versions of
$\ger G^0$ and $\ger G^1$, obtained by adjoining $2^{m^2}$ isolated 
vertices, not in the field of the pre-ordering $\preceq$.    
 
\begin{prop}		\label{cfi-cpt} 
There is a polynomial time BGS program that accepts $\ger H_m^0$ and 
rejects $\ger H_m^1$ for all $m$. 
\end{prop} 
 
\begin{pf} 
The program checks whether there is a choice of one vertex $(e,\pm)$ 
from each $U(e)$ such that each $U(v)$ contains a vertex whose 
neighbors were all chosen.  We saw above that such a choice is 
possible in $\ger H_m^0$ but not in $\ger H_m^1$.  To write the 
program in BGS, think of it as consisting of two phases.  In the first 
phase, it goes through all the blocks $U(e)$ in order (i.e., in the 
order induced on edges $e$ by $\preceq$), splitting into parallel 
subcomputations each of which has one choice of vertices from the 
$U(e)$'s.  Then in the second phase, each of these subcomputations 
goes through the $U(v)$'s in order, checking whether there is a vertex 
whose neighbors were all chosen. 
 
The number of edges in $G$ is $(m+1)m/2\leq m^2$.  So the number of 
parallel subcomputations is no bigger than $2^{m^2}$.  The padding in 
the definition of $\ger H_m^i$ ensures that this is a polynomial in 
the input size.  It easily follows that a PTime version of this BGS 
program does what is required of it in the proposition. 
\end{pf} 
 
In contrast, Lemma~2.8 of \cite{otto} (see also Corollary~7.1 of 
\cite{cfi}) shows that $\ger H_m^0$ and $\ger H_m^1$ cannot be 
distinguished by any sentence in $C^m_{\infty,\omega}$.  (Our graphs 
differ from Otto's in two ways.  Where we have a single vertex $(e,+)$ 
or $(e,-)$ joined to vertices in two $U(v)$'s corresponding to the two 
endpoints of $e$, he has two adjacent vertices, each joined to 
vertices in just one $U(v)$.  And his graphs are not padded with 
isolated vertices.  Neither of these differences affects the proof 
that the graphs are $C^m_{\infty,\omega}$-equivalent.) 
 
Thus, these examples of Cai, F\"urer, and Immerman, with sufficient 
padding, show that FP+Card does not include \cpt; so \cptc\ properly 
includes FP+Card. 
 
Since the padding looks very artificial, it is natural to ask what 
happens if we omit it.  Writing $\ger G_m^i$ for $\ger G^i$ in the 
special case where $G$ is a complete graph on $m+1$ vertices, we 
still have that $\ger G_m^0$ and $\ger G_m^1$ are 
$C^m_{\infty,\omega}$-equivalent, just as before.  But the proof of 
Proposition~\ref{cfi-cpt} breaks down, since the computation is no 
longer PTime bounded.  The input structures $\ger G_m^i$ have size 
only $(m+1)(2^{m-1}+m)$, so polynomial time would mean time bounded by 
$2^{cm}$ for some $c$.  This is insufficient for generating all the 
choices of $(e,\pm)$'s.  We do not know whether \cpt\ or even \cptc\ 
can distinguish all the $\ger G_m^0$'s from the $\ger G_m^1$'s, but an 
argument from the proof of Corollary~7.1 of \cite{cfi} can be adapted 
to give the following result. 
 
\begin{prop}		\label{cfi-ptime} 
The isomorphism closure of the class $\{\ger G_m^0:m\in\bbb N\}$ is in 
PTime. 
\end{prop} 
 
\begin{pf} 
We must exhibit a PTime algorithm which, given a structure \ger X of 
the appropriate vocabulary and given an ordering of its underlying 
set, decides whether $\ger X\cong\ger G_m^0$ for some $m$.    
 
It is straightforward to check whether $\ger X\cong\ger G_m^i$ for some
$m$ and some $i$: First count the number of vertices and use it to
compute
$m$.  Then check whether $\preceq$ is a linear pre-ordering with $m+1$
equivalence classes $U(v)$, each of size $2^{m-1}$.  Then check whether
the remaining vertices come in pairs $U(e)$, one for each pair of $v$'s.
Label the vertices in each pair $U(e)$ with $+$ and $-$, say using $+$
for
the earlier one in the given ordering of the set of vertices.  Then
label
each vertex in each $U(v)$ by the sequence of $+$'s and $-$'s describing
which vertices from the $U(e)$'s it is adjacent to.  (We can do this,
with
\emph{sequences} because of the ordering of vertices.)  Then check
whether
the sequences associated to any two distinct vertices from the same
$U(v)$
differ in a nonzero even number of locations.  \ger X has the form $\ger
X\cong\ger G_m^i$ if and only if all these computations and checks
succeed.
 
It remains to distinguish $\ger G_m^0$ from $\ger G_m^1$.  This is 
done by a slight variant of the approach used above in the padded case 
(Proposition~\ref{cfi-cpt}).  From each pair $(e,\pm)$ choose the 
vertex labeled $(e,-)$ above.  Label a block $U(v)$ 
``good'' if it contains a vertex adjacent only to chosen vertices and 
``bad'' otherwise.  Of course, if all blocks are good, then our graph 
\ger X is, up to isomorphism, $\ger G_m^0$.    
 
The same holds if the number of bad blocks is even.  Indeed, in this 
case, there is a set $S$ of edges of the complete graph $G$ such that 
\od S is exactly the set of vertices corresponding to bad blocks.  If 
we choose $(e,+)$ instead of $(e,-)$ at the edges $e\in S$, then the 
new choices have the property that every $U(v)$ contains a vertex 
adjacent only to chosen ones. 
 
On the other hand, if the number of bad blocks is odd, then a 
similar argument shows that $\ger X\cong\ger G_m^1$.    
 
It remains to observe that we can determine in PTime which blocks are 
bad and (thanks to the ordering) how many of them there are.  So this 
algorithm works in PTime and accepts precisely (the isomorphs of) the 
graphs $\ger G_m^0$. 
\end{pf} 
 
\begin{ques} 
Can \cptc\ (or even \cpt) distinguish $\ger G_m^0$ from $\ger G_m^1$ 
for all $m$? 
\end{ques} 
 
If the answer is negative, then we have a separation of PTime from 
\cptc.  If the answer is affirmative, then we merely have another 
separation of \cptc\ from FP+Card, which we already had using $\ger 
H_m^i$.  The new separation would be aesthetically preferable, since 
it avoids padding. 
 
\section{Multipedes}		\label{mult} 
 
In this section, we study Boolean queries concerning the multipedes 
introduced in \cite{multipede}.  We begin by recalling the relevant 
definitions and results from \cite{multipede}.  That paper uses five 
notions of $k$-multipede, for $k=1,2^-,2,3,4$.  We shall not need the 
first two of these, so we begin with 2-multipedes. 
 
\begin{df} 
A \emph{2-multipede} is a finite 2-sorted structure, the two sorts 
being called ``segments'' and ``feet,'' with the following data. 
\begin{ls} 
\item A function $S$ from feet to segments, such that every segment is 
the image of exactly two feet. 
\item A family of 3-element sets of segments, called ``hyperedges.'' 
(This family is coded as a totally irreflexive and symmetric ternary 
relation.) 
\item A family of 3-element sets of feet, called ``positive triples,'' 
(similarly coded). 
\end{ls} 
These data are subject to the following requirements. 
\begin{ls} 
\item If $P$ is a positive triple of feet, then its image $S(P)$ is 
a hyperedge.  In particular, $S$ is one-to-one on $P$. 
\item If $H$ is a hyperedge then, of the eight triples of feet that 
$S$ maps onto $H$, exactly four are positive. 
\item If $P$ and $P'$ are positive triples of feet with $S(P)=S(P')$, 
then $|P\bigtriangleup P'|$ is even. 
\end{ls} 
\end{df} 
 
Notice that, for any three-element set $H$ of segments, there are 
exactly eight three-element sets of feet mapped onto $H$ by $S$. 
These eight are partitioned into two sets of four by the equivalence 
relation ``even symmetric difference.''  The positivity relation picks 
out one of these two equivalence classes for each hyperedge $H$. 
 
\begin{df} 
A \emph{3-multipede} is a 2-multipede together with a linear ordering 
$\leq$ of the set of all segments. 
\end{df} 
 
\begin{df} 
A \emph{4-multipede} is a 3-multipede together with a third sort, 
called ``sets,'' and a binary relation $\eps$ between segments and 
sets such that every set (in the ordinary sense) of segments is 
$\{s:s\eps x\}$ for a unique set (in the sense of the structure) $x$. 
\end{df} 
 
In other words, up to isomorphism, the sort of sets is exactly the 
power set of the sort of segments and $\eps$ is the membership 
relation.   
 
\begin{df} 
A multipede is \emph{odd} if, for every nonempty set $X$ of segments, 
there is a hyperedge whose intersection with $X$ has odd cardinality. 
\end{df} 
 
The value of oddness and of the linear ordering in the definition of 
3-multipedes is the following result, combining Lemmas~4.1 and 4.4 of 
\cite{multipede}.   
 
\begin{prop}[\cite{multipede}] \label{rigid} 
Odd 3-multipedes and odd 4-multipedes are rigid. 
\end{prop} 
 
\begin{pf} 
Consider first an automorphism $\alpha$ of an odd 3-multipede. 
Because of the linear ordering, it must fix every segment.  So all it 
can do with feet is to interchange the two feet in $S^{-1}(\{s\})$ for 
certain segments $s$.  Let $X$ be the set of segments $s$ whose two 
feet $\alpha$ interchanges.  Since $\alpha$ preserves positivity of 
triples of feet, the intersection of $X$ with each hyperedge must have 
even cardinality.  Since the multipede is odd, this means $X=\emp$, 
and so $\alpha$ fixes all feet. 
 
In the case of a 4-multipede, we see as above that an automorphism 
fixes all segments and all feet.  In order to preserve $\eps$, it must 
also fix all sets. 
\end{pf} 
 
The main work in \cite{multipede} involves the notion of a $k$-meager 
multipede (for $k\in\bbb N$); we omit the definition here because we 
shall avoid needing it.  We do need two trivial (given the definition) 
and two deep properties of meagerness.  The trivial properties are 
that meagerness depends only on the segments and hypergraphs and that 
$k$-meagerness implies $l$-meagerness for all $l<k$.  The deep 
properties are the following two results from \cite{multipede}; the 
first is Theorem~3.1 and the second combines Lemmas~4.2 and 4.5 of 
\cite{multipede}.   
 
\begin{prop}[\cite{multipede}]	\label{meager} 
For any positive integers $l\geq2$ and $N$, there exists an odd,   
$l$-meager multipede with more than $N$ segments. 
\end{prop} 
 
We observe that it doesn't matter in this proposition whether 
``multipede'' refers to 2-, 3-, or 4-multipedes.  Once we have an odd, 
$l$-meager 2-multipede, we can expand it with an arbitrary linear 
ordering of its segments and we can adjoin an appropriate universe of 
sets to get an odd, $l$-meager 4-multipede. 
 
\begin{prop}[\cite{multipede}] 
No formula of $C^l_{\infty,\omega}$ can distinguish between the two 
feet of any segment in an $l$-meager multipede. 
\end{prop} 
 
The purpose of these results in \cite{multipede} was to exhibit a
finitely
axiomatizable (in first-order logic) class of structures, namely the odd
4-multipedes, such that all structures in the class are rigid but no
$C^\omega_{\infty,\omega}$ formula can define a linear ordering on all
structures of the class.  The addition of sets, in going from
3-multipedes
to 4-multipedes, served to make ``odd'' first-order definable.
 
In the present paper, our interest is in definability (or 
computability) of Boolean queries, not linear orderings.  To apply the 
ideas of \cite{multipede} in this context, we make one additional 
definition, intended to apply to 2-, 3,- and 4-multipedes 
simultaneously.   
 
\begin{df} 
A \emph{multipede with a shoe} is a multipede with a distinguished
foot, called the ``foot with a shoe'' or simply the ``shoe.'' In the
case of 3- and 4-multipedes, it is further required that $S$ of the
shoe is the first segment in the order $\leq$.
\end{df} 
 
The first Boolean query we shall consider is the isomorphism problem 
for 4-multipedes with shoes.  The input here is a pair of 
4-multipedes, each with a shoe.  (Since 4-multipedes are 3-sorted 
structures, it is convenient to regard a pair of them as a 6-sorted 
structure.)  The question is whether the two are isomorphic. 
 
\begin{thm}		\label{iso4easy} 
The isomorphism problem for 4-multipedes with shoes is in \cpt. 
\end{thm} 
 
\begin{pf} 
Since the segments of a 4-multipede are linearly ordered, any
isomorphism is uniquely determined on segments and therefore on sets,
and it is easy to check in choiceless polynomial time whether the
hyperedges in the two multipedes match up properly.  The only real
problem is whether the feet can be matched up so as to preserve $S$
and positivity.
 
If the input multipedes have $n$ segments each, then there are $2^n$ 
ways to match up the feet while preserving $S$, since for each of the 
$n$ pairs of corresponding segments in the two multipedes, there are 
two ways to match up their feet.  The problem is whether any of these 
$2^n$ matchings preserves positivity. 
 
Because of the universe of sets in a 4-multipede, the input structures 
are larger than $2^n$.  So a PTime bounded BGS algorithm has enough 
time to construct, in parallel, all the relevant matchings of feet and 
to check whether any of them preserve positivity. 
\end{pf} 
 
In this proof, the role of the sets is to serve as padding, making 
``polynomial time'' long enough to carry out the algorithm.  The only 
reason we didn't have to resort to explicit padding here (as we did in 
the case of the Cai, F\"urer, Immerman examples above) is that the 
necessary padding was already done, for a different purpose, in 
\cite{multipede}.  Of course, this raises the question of what happens 
without padding, i.e., with 3-multipedes; we shall return to this 
question after the next result, which completes our discussion of 
4-multipedes. 
 
\begin{thm}		\label{iso4hard} 
The isomorphism problem for 4-multipedes with shoes is not in 
$C^\omega_{\infty,\omega}$ and therefore not in FP+Card. 
\end{thm} 
 
\begin{pf} 
Suppose we had a sentence $\theta$ of $C^\omega_{\infty,\omega}$ 
expressing isomorphism between 4-multipedes with shoes.  Fix $l$ so 
large that $\theta$ is in $C^l_{\infty,\omega}$, and let $M$ be an 
odd, $l$-meager 4-multipede, which exists by Proposition~\ref{meager}. 
Let $\ger M_0$ and $\ger M_1$ be the two expansions of $M$ with shoes, 
i.e., one of the two feet of the first segment is the shoe in $\ger 
M_0$ and the other is the shoe in $\ger M_1$.  By 
Proposition~\ref{rigid}, $\ger M_0$ and $\ger M_1$ are not isomorphic, 
for an isomorphism would be a non-trivial automorphism of \ger M. 
Thus, $\theta$ must be false in the structure $\ger M_0+\ger M_1$ but 
true in $\ger M_0+\ger M_0$.    
 
This means that the Spoiler has a winning strategy in $C^l$ game for 
the pair of structures $\ger M_0+\ger M_1$ and $\ger M_0+\ger M_0$. 
(See \cite[Theorem~2.1]{otto}.)  We obtain a contradiction by 
exhibiting a winning strategy for the Duplicator in this game. 
 
The proofs of Lemmas~4.2 and 4.5 in \cite{multipede} provide a winning 
strategy for the Duplicator in the $C^l$ game for the pair of 
structures $\ger M_0$ and $\ger M_1$.  And the Duplicator has a 
trivial winning strategy for the pair $\ger M_0$ and $\ger M_0$: just 
copy whatever the Spoiler does.  Combining these two known strategies, 
we get a winning strategy for the Duplicator for the pair $\ger 
M_0+\ger M_1$ and $\ger M_0+\ger M_0$ as follows.  When Spoiler picks 
a subset of one of these structures, think of it as two subsets, one 
in each of the two component multipedes.  Apply the known strategies 
to find two subsets of the same cardinalities in the component 
multipedes of the other board, and play the union of these two 
subsets.  Then, when Spoiler picks a point in one of these subsets, 
pick a point on the other board by consulting the appropriate one of 
the known strategies.    
 
In effect, Duplicator is playing the $C^l$ game for $\ger M_0+\ger 
M_1$ and $\ger M_0+\ger M_0$ by playing separately the trivial game 
for the first components $\ger M_0$ and $\ger M_0$ and the game for 
the second components $\ger M_0$ and $\ger M_1$.  Since he wins in 
both components, he also wins the overall game. 
\end{pf} 
 
The last two theorems give us, once again, a \cpt\ computable query 
that goes beyond FP+Card. 
 
Turning to 3-multipedes, we see that Theorem~\ref{iso4hard} remains 
true with the same proof.  But the proof of Theorem~\ref{iso4easy} no 
longer applies, because a 3-multipede with $n$ segments has only $3n$ 
elements and so polynomial time is inadequate for producing all 
possible matchings of the feet.  As a result, we do not know whether 
isomorphism of 3-multipedes with shoes is in \cpt\ or even in \cptc. 
But we do have the following weaker result. 
 
\begin{thm}		\label{iso3} 
Isomorphism of 3-multipedes with shoes is computable in PTime. 
\end{thm} 
 
\begin{pf} 
We must present a PTime algorithm which, given a structure $\ger 
A+\ger B$ where \ger A and \ger B are 3-multipedes with shoes, and 
given a linear ordering $\preceq$ of their union, decides whether they 
are isomorphic.  There is a slight possibility of confusion between 
the different orderings here, the linear orderings of segments that 
are part of the 3-multipede structure of \ger A and \ger B, and the 
additional ordering $\preceq$ of the whole combined structure.  The 
latter will be used only to distinguish between the two feet of any 
segment (in either component multipede); we'll call one the left and 
the other the right foot.  We fix the terminology so that the left 
foot is $\prec$ the right \emph{except} that in both multipedes the 
shoe is declared to be the left foot of its segment regardless of what 
$\preceq$ does.  For the rest of the proof, any mention of an ordering 
refers to the orderings of segments that are part of the 3-multipede 
structure. 
 
As in the proof of Theorem~\ref{iso4easy}, thanks to the orderings of 
segments, there is no difficulty deciding whether the hypergraph 
structures on the segments agree.  The problem is to decide whether 
the feet can be matched appropriately.  If there are $n$ segments then 
there are $2^n$ possible (i.e., respecting the function $S$) 
matchings, and the algorithm lacks the time to check each one to see 
if it preserves positivity.  Instead, let the algorithm proceed as 
follows.   
 
Take one specific matching $\mu$, namely the one that maps left feet 
to left feet and (therefore) right feet to right feet (of the 
corresponding segments, of course).  If it happens to preserve 
positivity, then output ``yes'' (or ignore this obvious answer and 
proceed as in the general case).  For each of the two multipedes \ger 
A and \ger B, list all its hyperedges in lexicographic order (with 
respect to the ordering of segments).  Of course, the two multipedes 
have the same number, say $m$, of hyperedges; otherwise, the algorithm 
would have detected non-isomorphism earlier and we wouldn't be looking 
for a matching of the feet.  Form an $m$-component vector $\vec v$ of 
0's and 1's, where the $k^{\text{th}}$ entry is 0 if $\mu$ preserves 
positivity of triples of feet at the $k^{\text{th}}$ hyperedge and 1 
otherwise. 
 
Any other possible matching is obtained from $\mu$ by reversals at 
some set $X$ of segments of \ger A.  Call the result $\mu_X$.  (So 
$\mu=\mu_\emp$.)  In order for $\mu_X$ to be an isomorphism, i.e., to 
preserve positivity at all segments, $X$ must have an odd intersection 
with those hyperedges where $\mu$ failed to preserve positivity and an 
even intersection with the other hyperedges.  We reformulate this 
criterion as follows.  Represent any $X$ by an $n$-component vector of 
zeros and 1's, where the $k^{\text{th}}$ entry is 1 if and only if the 
$k^{\text{th}}$ segment is in $X$.  Also, let $A$ be the 
segment-hyperedge adjacency matrix; it is the $m\times n$ matrix whose 
$(k,l)$ component is 1 if the $k^{\text{th}}$ hyperedge contains the 
$l^{\text{th}}$ segment.  Then the condition for $\mu_X$ to preserve 
positivity is simply that $A\vec x=\vec v$, where both $\vec x$ and 
$\vec v$ are considered as column vectors, and where arithmetic is 
done modulo~2.    
 
Thus, the isomorphism question is reduced to the question of 
solvability of a system of linear equations $A\vec x=\vec v$ over the 
field $\bbb Z/2$.  But such questions are easily solved in polynomial 
time, by Gaussian elimination. 
\end{pf} 
 
We repeat the main question left open by the results in this section. 
 
\begin{ques} 
Is isomorphism of 3-multipedes with shoes computable in \cptc\ or in 
\cpt?   
\end{ques} 
 
A negative answer to the \cptc\ version of the question would separate 
PTime from \cptc.  A positive answer would only give yet another 
separation of \cptc\ from FP+Card (without any unpleasant padding). 
 
\section{Determinants}	\label{det} 
 
The use of linear algebra modulo 2 in the proof of Theorem~\ref{iso3} 
suggests that this topic or more generally linear algebra over finite 
fields may lead to interesting problems at or near the border between 
PTime and \cptc.  In this section, we consider problems of this sort, 
related to computing determinants or at least deciding whether a given 
matrix has zero determinant. 
 
\subsection{Matrices and Determinants} 
 
The method of Gaussian elimination, i.e., reducing a matrix to echelon 
form by row or column operations, computes determinants of $n\times n$ 
matrices in $O(n^3)$ arithmetical operations.  When the matrix entries 
come from a fixed finite field (or commutative ring), this observation 
shows that determinants are computable in polynomial time.  (If the 
matrix entries come from an infinite field or ring, then one must take 
into account how the entries are presented and how complex the 
arithmetical operations are.  We shall discuss the infinite case 
briefly below.) 
 
Matrices are usually regarded as having their rows and columns given in
a
specified order, and the Gaussian elimination algorithm makes use of
this
order in deciding which row operations to apply.  Our concern in this
section will be with ``matrices'' in which the rows and columns are
indexed by unordered sets; thus Gaussian elimination cannot be used.  We
use matrices as inputs to computations, so, in accordance with the
conventions of BGS, we shall code matrices as structures.
 
There are two inequivalent ways to make precise the notion of a matrix 
with unordered rows and columns.    
 
\begin{df}		\label{bipmatrix} 
Let $I$ and $J$ be finite sets, and let $R$ be a finite commutative 
ring.  An $I\times J$ \emph{matrix} with entries from $R$ is a 
function $M:I\times J\to R$.  We regard $M$ as a two-sorted structure, 
the sorts being $I$ and $J$, with basic relations 
$$ 
M_r=\{(i,j)\in I\times J:M(i,j)=r\} 
$$ 
for all $r\in R$. 
\end{df} 
 
\begin{df}		\label{sqmatrix} 
Let $I$ be a finite set, and let $R$ be a finite commutative ring.  An 
$I$-\emph{square matrix} with entries from $R$ is a function 
$M:I\times I\to R$.  We regard $M$ as a one-sorted structure with 
underlying set $I$ and with basic relations 
$$ 
M_r=\{(i,j)\in I\times I:M(i,j)=r\} 
$$ 
for all $r\in R$. 
\end{df} 

An alternative but equivalent way to code matrices as structures would
be
to include $R$ as an additional sort, to have the matrix itself as a
function $I\times J\to R$ or $I\times I\to R$, and to include in the
vocabulary names for all members of $R$.
 
Notice that, even when $|I|=|J|$, an $I\times J$ matrix differs in an
essential way from an $I$-square matrix.  An ordering of (the
structure representing) an $I\times J$ matrix independently orders
both $I$ and $J$; an ordering of (the structure representing) an
$I$-square matrix merely orders $I$.
 
Thus, an $I$-square matrix has a well-defined determinant in the 
following sense.  If one linearly orders $I$ then one obtains a matrix 
in the usual sense.  The determinant of this matrix is independent of 
the ordering because if one changes the ordering the effect is to 
permute the rows and the columns in the same way.  If the row 
permutation is odd and therefore reverses the sign of the determinant, 
then the column permutation reverses the sign again, restoring the 
original value.  In contrast, even when $|I|=|J|$, the determinant of 
an $I\times J$ matrix is defined only up to sign.  One gets a square 
matrix in the usual sense by fixing any orderings of $I$ and $J$, but 
changing to different orderings may change the sign of the 
determinant. 
 
We observe that the question whether a matrix has zero determinant 
makes good sense not only for $I$-square matrices but also for 
$I\times J$ matrices as long as $|I|=|J|$.  Although the determinant 
is defined only up to sign, the sign doesn't matter if we only care 
whether the determinant is zero.  Similarly, it makes good sense to 
speak of the rank of an $I\times J$ matrix (whether or not 
$|I|=|J|$).   
 
\begin{rmk}		\label{square-bij} 
One can view an $I$-square matrix as an $I\times J$ matrix together 
with a specified bijection between $I$ and $J$.  Every structure of 
the latter sort is isomorphic to one where $I=J$ and the specified 
bijection is the identity; if the isomorphism is required to be the 
identity on $I$ then it is unique. 
\end{rmk} 
 
\subsection{Determinants Modulo Two} 
 
In this subsection, we consider determinants of square matrices with 
entries from the two-element field $\bbb Z/2$.  For this particular 
field, an $I$-square matrix $M$ can be regarded as a directed graph 
with vertex set $I$ and arc set $M_1=\{(i,j)\in I^2:M(i,j)=1\}$, for 
the other relation, $M_0$, in the structure representing $M$ is then 
determined as the complement of $M_1$.  In other words, any square 
matrix over $\bbb Z/2$ can be regarded as the incidence matrix of a 
directed graph.  The graph here may have loops and may have pairs of 
opposite arcs $(i,j)$ and $(j,i)$ but cannot have parallel arcs; an 
arc is simply an ordered pair of vertices. 
 
Another simplification resulting from the restriction to $\bbb Z/2$ is 
that the problems ``compute the determinant'' and ``is the determinant 
zero?'' are equivalent, since there is only one possible non-zero 
value.  We shall consider the problem in the form ``is the determinant 
zero,'' for it is in this form that our results generalize to other 
finite fields.    
 
A third simplification is that determinants are well-defined for 
$I\times J$ matrices with $|I|=|J|$.  The sign ambiguity described 
earlier disappears in characteristic 2 where $x=-x$.  Nevertheless, 
the algorithm presented in this subsection applies to $I$-square 
matrices only.  From the point of view described in 
Remark~\ref{square-bij}, we shall make real use of the given bijection 
between the rows and the columns.  Later, we shall consider ways to 
avoid this. 
 
\begin{thm}		\label{det2+} 
The determinant of the square matrix over $\bbb Z/2$ represented by a 
finite directed graph is definable in FP+Card. 
\end{thm} 
 
Here we identify the possible values 0 and 1 of the determinant with 
the truth values, so that the determinant becomes a Boolean query. 
 
\begin{pf} 
We describe an algorithm for deciding whether any $I$-square matrix is 
non-singular.  The algorithm is easily seen to be formalizable as a 
polynomial time algorithm in BGS+Card.  Afterward we sketch how to 
convert the algorithm into a definition in FP+Card. 
 
We begin with a preliminary observation.  Given two $I$-square 
matrices $M$ and $N$, we can compute the product matrix $MN$, which is 
also an $I$-square matrix.  Indeed, $(i,j)$ is an arc in the graph 
$MN$ if and only if the cardinality of the set   
$$ 
\{k\in I: (i,k)\in M_1\text{ and }(k,j)\in N_1\} 
$$ 
is odd.  Since the parity of a natural number (which may be regarded
as a von~Neumann ordinal --- see Section~\ref{background}) is easily
in \cpt, it follows that all entries of the product matrix can be
computed in \cptc.
 
Next, we observe that we can compute powers of a matrix, even when the
exponent is so large that it is given in binary notation.  We first
describe how binary notation for natural numbers can be handled in the
BGS context.

The idea is that the binary representation of a natural number $r$, say
of
length $l=\text{lg}(r)$, amounts to a subset $C$ of $\{0,\dots,l-1\}$,
namely the set of places where a 1 occurs in the binary notation.  Thus
$r=\sum_{c\in C}2^c$.  We remark that, for non-zero $r$ and therefore
nonempty $C$, the largest element of $C$ is easily computable from $C$,
namely as $\bigcup C$ (where natural numbers are identified with
von~Neumann ordinals).

Suppose we are given an $I$-square matrix $M$ and an integer $r$ in
binary notation, i.e., the set $C$ as above.  Then we can compute
$M^r$ in time polynomial in $|I|$ and $\text{lg}(r)$.  The
computation of $M^r$ is done by repeated squaring, i.e., by applying
the recursion formulas
$$ 
M^r= 
\begin{cases} 
M&\text{if }r=1\\ 
\left(M^{r/2}\right)^2&\text{if }r\geq2\text{ is even}\\ 
\left(M^{(r-1)/2}\right)^2\cdot M&\text{if }r\geq2\text{ is odd.} 
\end{cases} 
$$ 
Here is  a BGS program for this algorithm, using matrix multiplication 
as an ``external'' function, which means that for the complete 
algorithm one should replace all matrix multiplications here by the 
algorithm described above.
 
\begin{eatab} 
do in parallel\\ 
\>if $Mode=0$ then \\ 
\>\>do in parallel $X:=M$; $p=1+\max C$; $Mode:=1$ enddo\\ 
\>endif;\\ 
\>if $Mode=1$ and $p=0$ then Halt:=true endif\\ 
\>if $Mode=1$ and $p\neq0$ and $p-1\in C$ then\\ 
\>\>do in parallel $X:=X\cdot X\cdot M$; $p:=p-1$ enddo\\ 
\>endif\\ 
\>if $Mode=1$ and $p\neq 0$ and $p-1\notin C$ then\\ 
\>\>do in parallel $X:=X\cdot X$; $p:=p-1$ enddo\\ 
\>endif\\ 
enddo 
\end{eatab} 
 
As a final preparatory step, we compute the order of $\text{GL}_n(\bbb 
Z/2)$, the group of non-singular $n\times n$ matrices over $\bbb 
Z/2$.  This order is   
$$ 
g=(2^n-1)(2^n-2)(2^n-4)\cdots(2^n-2^{n-1}) 
=\prod_{i=0}^{n-1}(2^n-2^i). 
$$ 
To see this, we use the fact that an $n\times n$ matrix is 
non-singular if and only if its columns are linearly independent 
vectors in $(\bbb Z/2)^n$.  If we imagine the columns being chosen one 
at a time, the first column of such a matrix can be any non-zero 
vector in $(\bbb Z/2)^n$; the second can be any vector different from 
the first and from 0; the third can be any vector that is not a linear 
combination of the first two; and in general any column can be any 
vector not a linear combination of the previously chosen columns. 
Thus, there are $2^n-1$ choices for the first column, each leaving 
$2^n-2$ choices for the second, each leaving $2^n-4$ choices for the 
third, and so on. 
 
Now given an $I$-square matrix $M$, in the form of a digraph with 
vertex set $I$, we can determine whether it is non-singular as 
follows in \cptc.  First, use the cardinality function to determine 
the von~Neumann ordinal $n=|I|$.  From this, compute the group order 
$g$ in binary notation.  Notice that $g<2^{n^2}$, so the length of 
this binary expansion ($\max C$ in the notation above) is bounded by 
$n^2$.  Our formula for $g$ above makes the computation of this binary 
expansion a simple matter, easily programmed in BGS (without further 
use of the cardinality function).  Next, compute $M^g$; as indicated 
above, this can be done in \cptc.  Finally, output 1 if $M^g$ is the 
identity matrix (i.e., if the arcs in the digraph $M^g$ are exactly 
the loops $(i,i)$ for all $i\in I$) and 0 otherwise. 
 
To see that this algorithm gives the correct answer, recall from 
elementary group theory the fact (a special case of Lagrange's 
theorem) that the order of an element in a group always divides the 
order of the group.  Thus, if $M$ is non-singular then the matrix 
obtained by ordering $I$ arbitrarily is an element of 
$\text{GL}_n(\bbb Z/2)$, so its $g^{\text{th}}$ power is the identity 
matrix, and the same follows for $M$.  If, on the other hand, $M$ is 
singular, then so are all its powers; in particular none of its powers 
is the identity matrix.  This completes the proof that non-singularity 
of square matrices over $\bbb Z/2$ is computable in \cptc. 
 
Finally, we briefly indicate why this algorithm yields a definition in 
FP+Card.  Since the input is a structure (directed graph) of size $n$, 
FP+Card works with a two-sorted structure $\ger A^*$ consisting of the 
input graph and the natural numbers up to $n$.  The algorithm above 
used natural numbers up to $n^2$ (to produce the binary expansion of 
$g$), but these can be coded as pairs of numbers below $n$.  The 
computation of $g$ (in binary form) is a polynomial time algorithm 
working on a numerical input (the second sort of $\ger A^*$), so it 
can be expressed in FP.  The repeated squaring algorithm for computing 
$M^g$ can be cast as a definition, using the fixed-point operator, of 
the ternary relation   
$$ 
\{(i,j,q):(i,j)\in M^{g_q}\} 
$$ 
where $g_q$ means the integer represented by the $q$ most significant 
digits in the binary expansion of $g$.  (More precisely, this is a 
quaternary relation because, as indicated above, $q$ is represented by 
a pair of elements of the numerical sort in $\ger A^*$.)  Finally, the 
comparison between $M^g$ and the identity matrix is expressible in 
first-order logic. 
\end{pf} 
 
To complement the previous theorem, we show next that the cardinality 
function is essential in this or any choiceless algorithm for 
determinants over the two-element field. 
 
\begin{thm}		\label{det2-} 
The determinant of the square matrix over $\bbb Z/2$ represented by a 
finite directed graph is not computable in \cpt. 
\end{thm} 
 
\begin{pf} 
Temporarily fix a positive integer $n$.  As we saw in the proof of 
Theorem~\ref{det2+}, the number of non-singular $n\times n$ matrices 
over $\bbb Z/2$ is   
$$ 
g=\prod_{i=0}^{n-1}(2^n-2^i). 
$$ 
Since the total number of $n\times n$ matrices over $\bbb Z/2$ is 
$2^{n^2}$, the probability that such a matrix, chosen uniformly at 
random, is non-singular is   
$$ 
\frac g{2^{n^2}}=\prod_{i=0}^{n-1}\left(1-\frac{2^i}{2^n}\right) 
=\prod_{j=1}^{n}\left(1-\frac1{2^j}\right). 
$$ 
This product is therefore the probability that a random (with respect 
to the uniform distribution) directed graph on an $n$-element vertex 
set has, when viewed as a matrix, determinant 1. 
 
Now un-fix $n$ and let it tend to infinity.  The asymptotic 
probability that a large, random, directed graph has determinant 1 is 
$$ 
\prod_{j=1}^{\infty}\left(1-\frac1{2^j}\right). 
$$ 
This infinite product is obviously strictly smaller than 1.  It is
strictly greater than 0 (i.e., it converges in the conventional
terminology) because the series $\sum_j(1/2^j)$ converges.  (Recall
the standard proof: $1-x>e^{-2x}$ for all positive
$x\leq\frac12$.  Apply this to $x=1/2^j$ and take the product over
$j$, obtaining a convergent sum in the exponent.)
 
But the zero-one law proved by Shelah \cite{sh634} (see also 
\cite{01law}) implies that any property of digraphs computable in 
\cpt\ must have asymptotic probability 0 or 1.  Therefore, 
``determinant 1'' is not such a property. 
\end{pf} 
 
\subsection{Other Finite Fields} 
 
The FP+Card definition of ``non-singular'' given in the preceding 
subsection for square matrices over $\bbb Z/2$ works, with minor 
modifications, over any finite field $F$.  Of course, when the field 
has more than two elements, the ``non-singular'' question is weaker 
than the problem of actually evaluating the determinant. 
 
To indicate the minor modifications explicitly, let $F$ be a finite 
field of characteristic $p$ and cardinality $q=p^e$.  Then to decide 
non-singularity of square matrices over $F$, we can use the algorithm 
described above for the special case $q=2$ with the following two 
changes.  First, the order of the group $\text{GL}_n(F)$ is 
$$ 
g=(q^n-1)(q^n-q)(q^n-q^2)\dots(q^n-q^{n-1}) 
=\prod_{i=0}^{n-1}(q^n-q^i); 
$$ 
i.e., $q$ replaces 2 in the earlier formula.    
 
Second, multiplying matrices becomes slightly more tedious but remains 
straightforward.  Given two $I$-square matrices $M$ and $N$, to 
compute the $(i,j)$ entry of a product matrix $MN$, first do the 
following for each element $z\in F$.  Consider the set $P_z\subseteq 
F^2$ of pairs $(x,y)$ whose product in $F$ is $z$.  Since $F$ is fixed 
in this discussion, our BGS program or FP+Card formula can contain a 
complete listing of all the $P_z$'s.  Use the cardinality function to 
obtain the numbers 
$$ 
m_z=|\{k\in I:(M(i,k),N(k,j))\in P_z\}| 
$$ 
and then, in a trivial polynomial time computation, reduce these 
numbers modulo $p$ to obtain $\bar m_z=m_z\bmod p$.  Then the $(i,j)$ 
entry of $MN$ is the element of $F$ given by the sum 
$$ 
\sum_{z\in F}m_z\cdot z=\sum_{z\in F}\bar m_z\cdot z. 
$$ 
Since there are only finitely many ($p^q$) possible functions 
$z\mapsto \bar m_z$, our program or formula can contain a table 
giving, for each of these functions, the value of the sum.    
 
The following proposition summarizes the preceding discussion. 
 
\begin{prop}		\label{detp+} 
For any finite field $F$, there is an FP+Card definition of 
non-singularity for square matrices over $F$. 
\end{prop} 
 
At two points in the preceding discussion, we used that the field $F$ 
is fixed, so that our FP+Card formula can contain complete 
descriptions of the sets $P_z$ and the sums associated to the 
functions $z\mapsto\bar m_z$.  It is not difficult, however, to adjust 
the algorithm to work uniformly over all finite fields $F$, in time 
polynomial in $|F|=q$ and the size of the matrix.  In the first place, 
the table of all the $P_z$'s is essentially the multiplication table 
of the field; its size is only quadratic in $q$. So this table can be 
computed as part of the algorithm. 
 
There isn't enough time to compute the sums associated to all possible 
functions $z\mapsto\bar m_z$, since there are $p^q$ of these 
functions.  But when multiplying a particular pair of matrices, we 
need the sum for only one such function per entry.  Each single sum is 
easy to compute provided we are given an ordering of $F$.  So there is 
no difficulty computing, in polynomial time, the $n^2$ sums actually 
needed.  Thus, we obtain the following uniform version of the 
preceding proposition. 
 
\begin{prop}		\label{finf} 
There is an FP+Card formula defining non-singularity of square 
matrices over finite fields, where the input structure consists of a 
finite field $F$, a linear ordering of the set $F$, a set $I$, and an 
$I$-square matrix $M:I^2\to F$. 
\end{prop} 
 
To avoid possible confusion, we point out that there is no necessary 
connection between the linear ordering of $F$ and the field 
operations. 
 
\begin{ques} 
Can the determinant of a square matrix over a finite field be computed 
in \cptc?  Can it be defined in FP+Card?    
\end{ques} 
 
\subsection{Integer Matrices} 	\label{intmat}
 
In this subsection, we apply the preceding results to matrices with
entries from the ring \bbb Z of integers.  Since we require inputs of
computations to be finite structures for finite vocabularies, we must
modify the representation of matrices as structures described in
Definitions~\ref{bipmatrix} and \ref{sqmatrix}.  Those definitions
would yield an infinite vocabulary whenever the underlying ring is
infinite, and if we represented the matrix by a function (as in
Proposition~\ref{finf}) instead of a family of relations then the
vocabulary would be finite but the underlying set of the structure
would be infinite.  We adopt the convention that matrix entries are to
be written in binary notation.  Recall that this means that an entry
$r$ is represented by a set $C$ of natural numbers, the set of
locations of ones in the binary expansion.  Thus, each matrix entry is
to be a set of natural numbers, and therefore the matrix itself
amounts to a ternary relation, $M(i,j,s)$ with the meaning ``the
coefficient of $2^s$ in the binary expansion of the $(i,j)$ entry of
$M$ is 1.'' There are two problems with this set-up.

The smaller problem is that we have not taken into account the signs
of the matrix entries.  So we shall need a second relation, a binary
one, with the meaning ``the $(i,j)$ entry of $M$ is positive.''

The more serious problem is that, although the first and second
arguments
of the ternary relation $M$ are atoms, namely indices for rows or
columns
of our matrix, the third argument is a natural number.  Both BGS and
FP+Card are set up so that the numbers (von~Neumann ordinals in the case
of BGS, the numerical second sort in the case of FP+Card) are not part
of
the input structure.  So $M$ is not appropriate as an input in BGS or
FP+Card.  We therefore include in the input structure a copy of enough
of
the natural number system to allow coding our binary numbers.  Our
official representation of integer matrices will thus involve surrogate
natural numbers $\hat0,\hat1,\dots,\hat k$, although for practical
purposes, it does no harm to think of $0,1,\dots, k$ instead.

That is, an $I$-square matrix $M$ will be regarded as a two-sorted
structure with underlying sets $I$ and a set of indices
$\{\hat0,\hat1,\dots,\hat k\}$; the relations on this structure are
the linear ordering $\hat0<\hat1<\dots<\hat k$ on the second sort, the
ternary relation $M(i,j,\hat s)$ defined by ``the coefficient of $2^s$
in the binary expansion of the absolute value of the $(i,j)$ entry of
$M$ is 1,'' and the binary relation ``the $(i,j)$ entry of $M$ is
positive.''
 
The number $k$ in this representation of a matrix $M$ would ordinarily 
be taken as small as possible, so it is essentially the logarithm of 
the largest absolute value of the matrix entries.    
 
\begin{thm} 
There is an FP+Card formula which, on matrices $M$ over \bbb Z 
represented as structures as above, defines ``$M$ is non-singular.'' 
\end{thm} 
 
\begin{pf} 
We describe an algorithm for deciding whether a square matrix over 
\bbb Z is non-singular, and we show that it works, without arbitrary 
choices, in polynomial time.  The details of formalization in \cptc\ 
or in FP+Card will, however, be left to the reader. 
 
Given a square matrix $M$, represented as above by a structure with 
underlying sets $I$ (indexing the rows and columns) and 
$\{\hat0,\hat1,\dots,\hat r\}$ (indexing the digits in each entry), 
the algorithm proceeds as follows.  First, find the cardinalities 
$|I|$ and $r+1$ of these two sets, and let $n$ be the larger of the 
two.  Thus, the size of the matrix is at most $n\times n$ and each 
entry is at most $2^n$ in absolute value.  The entire algorithm will 
take time polynomial in $n$. 
 
Second, generate a list of the first $2n^2$ prime numbers.  (For a BGS
algorithm, the primes are represented as von~Neumann ordinals.  For an
FP+Card definition, they are represented by rather short tuples of
elements from the numerical sort; in fact, triples will suffice --- see
below.)  This list can be produced by applying the sieve of
Eratosthenes.
The time required by the sieve of Eratosthenes is polynomial relative to
the prime numbers involved (though it is not polynomial relative to the
lengths of the primes in binary notation).  And the primes involved here
are, according to the prime number theorem, below $n^3$ provided $n$ is
large enough.  Therefore the time needed to generate this list of primes
is polynomial in $n$.
 
Third, go through all the primes $p$ in the list, checking for each 
one whether $M$ reduced modulo $p$ is non-singular as a matrix over 
$\bbb Z/p$.  The results of the previous subsection show that this can 
be done in polynomial time. 
 
Finally, output ``yes,'' meaning that $M$ is non-singular, if and only 
if it was non-singular modulo at least one of the primes $p$ on the 
list.   
 
This algorithm can clearly be programmed in BGS with the cardinality 
function, and it runs in polynomial time.  It is a routine matter to 
formalize it in FP+Card.  It remains to show that it gives the correct 
answer.    
 
If $M$ is singular, i.e., if its determinant is zero, then its 
reduction modulo $p$ has determinant zero for every prime $p$. 
Conversely, suppose $M$ is non-singular, and let its determinant be 
$d\neq0$.  Of course, then $M$ is non-singular modulo some primes, for 
example any primes larger than $|d|$, but we must show that it is 
non-singular modulo at least one of the first $2n^2$ primes. 
 
For this purpose, we first estimate how big $|d|$ might be.  The 
determinant of a $k\times k$ matrix is the sum of $k!$ terms, each of 
which is the product of $k$ of the matrix's entries.  In our 
situation, this means that $d$ is the sum of at most $n!$ terms, each 
the product of at most $n$ numbers, each at most $2^n$ in absolute 
value.  Thus, each of these $n!$ products is at most $2^{n^2}$. 
Therefore, 
$$ 
|d|\leq n!\cdot 2^{n^2}\leq n^n\cdot 2^{n^2}\leq (2^n)^n\cdot 
 2^{n^2}=2^{2n^2}. 
$$ 
 
Recall that the list of primes used by our algorithm consisted of the 
first $2n^2$ primes.  Since each prime is $\geq2$, the product of the 
listed primes is larger than $2^{2n^2}$ and therefore larger than 
$|d|$.  Thus, $|d|$, being non-zero, cannot be divisible by this 
product of primes.  That means that at least one prime $p$ on our list 
fails to divide $|d|$, i.e., that $M$ is non-singular modulo $p$, and 
therefore the algorithm gives the correct answer. 
\end{pf} 
 
\begin{rmk} 
The algorithm described in the preceding proof does more than the 
theorem claims.  It determines exactly which primes divide the 
determinant $d$ of $M$.  Indeed, it checks this divisibility directly 
for the first $2n^2$ primes, and the proof shows that no larger prime 
can divide $d$ unless $d=0$. 
 
The algorithm does not quite determine the value of $d$, for it does 
not determine its sign nor does it determine, for primes $p$ dividing 
$d$, whether $p^2$ or higher powers also divide $d$.  We do not know 
whether the determinant of an integer matrix can be computed (as a 
signed binary expansion) in polynomial time by a BGS program with the 
cardinality function. 
\end{rmk} 
 
\subsection{Rows and Columns May Differ} 
 
In the preceding discussion of determinants, we have dealt only with 
$I$-square matrices.  Up to a sign, determinants make sense for 
$I\times J$ matrices when $|I|=|J|$, and it makes sense to ask whether 
non-singularity of such matrices can be computed in \cptc\ or defined 
in FP+Card.  The algorithms from the preceding subsections do not 
suffice for this purpose, for they depend on taking powers of the 
given matrix, and $M^2$ is well-defined only when the rows and columns 
of $M$ are indexed by the same set.  Nevertheless, these algorithms 
can be modified to work when the rows and columns are indexed by 
different sets of the same size. 
 
\begin{thm} 
There is an FP+Card formula defining non-singularity of matrices over 
finite fields, where the input structure consists of a finite field 
$F$, a linear ordering of the set $F$, two sets $I$ and $J$ with 
$|I|=|J|$, and an $I\times J$ matrix $M:I\times J\to F$. 
\end{thm} 
 
\begin{pf} 
Let $M$ be an $I\times J$ matrix as in the statement of the theorem. 
Although $M^2$ is not defined when $I\neq J$, $M\cdot M^t$ is defined, 
where the superscript $t$ means transpose.  Furthermore, $M\cdot M^t$ 
is an $I$-square matrix.  Its entry in position $(i,i')$ is 
$\sum_{j\in J} M(i,j)M(i',j)$ which makes good sense for any $i,i'\in 
I$.  So by Proposition~\ref{finf} we can define non-singularity of 
$M\cdot M^t$ by an FP+Card formula.  But this is the same as defining 
non-singularity of $M$, since the determinant of $M\cdot M^t$ is the 
square of the determinant of $M$. 
\end{pf} 
 
\begin{rmk} 
An alternative proof of the theorem uses, instead of $M\cdot M^t$, the 
block matrix 
$$ 
\begin{pmatrix} 
0&M\\M^t&0 
\end{pmatrix}. 
$$ 
If $I$ and $J$ are disjoint (otherwise replace them by $I\times\{0\}$ 
and $J\times\{1\}$) then this block matrix is an $(I\cup J)$-square 
matrix.  Its $(x,y)$ entry is 0, if $x$ and $y$ are both in $I$ or both 
in $J$; $M(x,y)$, if $x\in I$ and $y\in J$; and $M(y,x)$ if $x\in J$ 
and $y\in I$.  Proposition~\ref{finf} allows us to define 
non-singularity of this block matrix, but again this is the same as 
non-singularity of $M$, since the determinant of the block matrix is 
the square of the determinant of $M$. 
\end{rmk} 
 
\section{Open Problems} 
 
The main problem that remains open is whether there is a logic, in the 
sense of \cite{logic}, that captures polynomial time on unordered 
structures.  It was conjectured in \cite{logic} that the answer is 
negative. 
 
A special case of the main problem is whether \cptc\ captures PTime. 
Of course a negative answer here is even more likely, but we have not 
been able to prove it.  In view of the results in this paper, a 
negative answer for the special case would follow from a negative 
answer to any of the following questions. 
\begin{ls} 
\item Can a \cptc\ program distinguish between the (unpadded) Cai, 
F\"urer, Immerman graphs $\ger G_m^0$ and $\ger G_m^1$ (as defined in 
Section~\ref{cfi-graphs}) for all $m$? 
\item Can isomorphism of 3-multipedes with shoes be decided by a 
\cptc\ program? 
\item Can a \cptc\ program decide whether a given graph (not
necessarily bipartite) admits a complete matching?
\item Can a \cptc\ program compute, up to sign, the determinant of 
an $I\times J$ matrix over a finite field (where $|I|=|J|$)? 
\end{ls} 
We point out that, although we have formulated these questions for 
\cptc, the logic in which we are primarily interested, the first two 
of them are open also for \cpt, and the last two are open also for 
FP+Card.

\end{document}